\documentclass{amsart}
\usepackage{amsmath,amsthm}
\usepackage{amsfonts,amssymb}

\newtheorem{thm}{Theorem}[section]
\newtheorem{cor}[thm]{Corollary}
\newtheorem{lem}[thm]{Lemma}
\newtheorem{prop}[thm]{Proposition}

\newtheorem{defn}[thm]{Definition}

\theoremstyle{remark}
\newtheorem{rem}{Remark}[section]

 \def\CD{{\mathcal D}}
 
 \def\CH{{\mathcal H}}

 \def\CP{{\mathcal P}}

 \def\CV{{\mathcal V}}
 \def\CW{{\mathcal W}}

 \def\NN{{\mathbb N}}

 \def\RR{{\mathbb R}}
 \def\ZZ{{\mathbb Z}}
        
        \def\proj{\operatorname{proj}}

\begin{document}
 
\title
  {Weighted Approximation of functions on the unit sphere}
\author{Yuan Xu}
\address{Department of Mathematics\\ University of Oregon\\
    Eugene, Oregon 97403-1222.}\email{yuan@math.uoregon.edu}

\date{July 3, 2002, revised December 9}
\thanks{Work supported in part by the National Science 
Foundation under Grant DMS-0201669}
                                                    
\begin{abstract}
The direct and inverse theorems are established for the best approximation 
in the weighted $L^p$ space on the unit sphere of $\RR^{d+1}$, in which the 
weight functions are invariant under finite reflection groups. The theorems 
are stated using a modulus of smoothness of higher order, which is proved to 
be equivalent to a $K$-functional defined using the power of the spherical 
$h$-Laplacian. Furthermore, similar results are also established for weighted 
approximation on the unit ball and on the simplex of $\RR^d$.
\end{abstract}

\maketitle                      
 
\section{Introduction} 
\setcounter{equation}{0}

Let $S^d = \{x: \|x\| =1\}$ denote the unit sphere in $\RR^{d+1}$, where 
$\|x\|$ denote the usual Euclidean norm. In the literature, the best $L^p$ 
approximation by polynomials on $S^d$ usually deals with the norm defined 
with respect to the Lebesgue measure, which is the unique measure on the 
sphere invariant under the rotation group, and the spherical harmonics play 
an essential role in the study. In this paper we study the weighted best 
$L^p$ approximation for a family of weight functions that are invariant under
reflection groups. 

For a nonzero vector $v \in \RR^{d+1}$, let $\sigma_v$ denote the reflection 
with respect to the hyperplane perpendicular to $v$, $x \sigma_v : = x - 
2 (\langle x,v \rangle /\|v\|^2) v$, $x \in \RR^{d+1}$, where $\langle x,y
\rangle$ denote the usual Euclidean inner product. Let $G$ be a finite 
reflection group on $\RR^{d+1}$ with a fixed positive root system $R_+$, 
normalized so that $\langle v, v \rangle =2$ for all $v \in R_+$. Then 
$G$ is a subgroup of the orthogonal group generated by the reflections 
$\{\sigma_v: v \in R_+\}$. Let $\kappa$ be a nonnegative multiplicity function 
$v \mapsto \kappa_v$ defined on $R_+$ with the property that $\kappa_u = 
\kappa_v$ whenever $\sigma_u$ is conjugate to $\sigma_v$ in $G$; then 
$v \mapsto \kappa_v$ is a $G$-invariant function. We consider the weighted 
$L^p$ best approximation with respect to the measure $h_\kappa^2 d \omega$ 
on $S^{d}$, where $h_\kappa$ is defined by 
\begin{equation}\label{eq:1.1}
h_\kappa(x) = \prod_{v \in R_+} |\langle x, v\rangle|^{\kappa_v}, \qquad 
   x \in \RR^{d+1}, 
\end{equation}
and $d \omega$ is the surface (Lebesgue) measure on $S^d$. 
The function $h_\kappa$ is a positive homogeneous function of degree 
$\gamma_\kappa:= \sum_{v \in R_+} \kappa_v$, and it is invariant under the 
reflection group $G$. The simplest example is given by the case $G=
\ZZ_2^{d+1}$ for which $h_\kappa$ is just the product weight function
\begin{equation}\label{eq:1.2}
  h_\kappa (x)  = \prod_{i=1}^{d+1} |x_i|^{\kappa_i}, \qquad \kappa_i \ge 0. 
\end{equation}    
We denote by $a_\kappa$ the normalization constant of $h_\kappa$, 
$a_\kappa^{-1}  = \int_{S^d} h_\kappa^2(y) d\omega$, and denote by 
$L^p(h_\kappa^2)$, $1\le p\le \infty$, the space of functions defined on 
$S^d$ with the finite norm
$$
\|f\|_{\kappa,p} := \Big(a_\kappa \int_{S^d} |f(y)|^p h_\kappa^2(y) d\omega(y) 
\Big)^{1/p}, \qquad  1 \le p < \infty,
$$
and for $p = \infty$ we assume that $L^\infty$ is replaced by $C(S^{d})$, the 
space of continuous functions on $S^d$ with the usual uniform norm 
$\|f\|_\infty$. 

The homogeneous polynomials that are orthogonal with respect to $h_\kappa^2 
d\omega$ are called $h$-harmonics, they are defined and studied by Dunkl 
(\cite{D1,D2}; see \cite{DX} and the references therein). The $h$-harmonics
satisfy many properties that are similar to those of ordinary harmonics. 
In particular, results on summability of the $h$-harmonic expansions have 
been developed in \cite{X97b,LX,X01a,X02}. While some of the results can be 
derived using methods similar to those used for ordinary harmonics, others 
become much more difficult to establish, largely due to the fact that the 
orthogonal group acts transitively on the sphere $S^d$ but a reflection group
does not. 

In order to understand the situation, to detect the obstacle in extending 
results for the Lebesgue measure to measures invariant under the reflection 
group, we study best approximation in $L^p(h_\kappa^2)$ in this paper; that is,
we consider 
$$
   E_n(f)_{\kappa,p} : = \inf \left \{ \|f - P \|_{\kappa,p}: 
        P \in \Pi_n^{d+1} \right\}, 
$$
where $\Pi_n^{d+1}$ denote the space of polynomials of degree at most $n$ in 
$d+1$ variables. For the Lebesgue measure on $S^d$, the problem of best 
approximation has been studied by many authors. We refer to 
\cite{BBP,Kam,LN,P,Rus} 
and the references therein. Much of our study uses the ideas of these authors,
since it turns out that their results can be extended to the weighted case 
with a proper definition of the modulus of smoothness. Such a definition is 
given recently in \cite{X02} in terms of a weighted spherical means for
$h_\kappa^2$ in \eqref{eq:1.2}, where we 
also proved that the modulus of smoothness is equivalent to a K-functional. If 
fact, the starting point in \cite{X02} is the K-functional since it arises 
naturally from the study of the de la Vall\'ee Poussin means of the 
$h$-harmonic series (see also \cite{BL}). It turns out that this modulus of 
smoothness can be used to give both direct and inverse theorems for the best 
approximation in $L^p(h_\kappa^2)$. Furthermore, 
it is possible to extend its definition to modulus of smoothness of higher 
order, prove that the extension is equivalent to the K-functional of higher 
order, and use it to establish direct and inverse theorems. This extends
the full strength of the work for the Lebesgue measure on the sphere in 
Rustamov \cite{Rus} to the weighted case. 

As one consequence of a weighted approximation theory for the unit sphere
$S^d$, we are able to establish a similar theory for the weighted approximation
on the unit ball $B^d = \{x: \|x\|\le 1\}$ of $\RR^d$, in which the weight 
function takes the form
$$
 W_{\kappa,\mu}^B(x) = h_\kappa^2(x) (1-\|x\|^2)^{\mu-1/2}, \qquad x \in B^d,
$$
where $h_\kappa$ is a reflection invariant weight function on $\RR^d$ and
$\mu \ge 0$. The case $h_\kappa (x) =1$ corresponds to the classical weight
function $W_\mu(x) = (1-\|x\|^2)^{\mu-1/2}$. Moreover, there is also a close
relation between the unit ball and the simplex $T^d = \{x \in \RR^d:
x_1 \ge 0, \ldots, x_d \ge 0, 1-|x| \ge 0\}$, where $|x| = x_1+\ldots + x_d$,
which allows us to further extend the theory to the weight approximation on 
$T^d$, in which the weight functions take the form
$$
W_{\kappa,\mu}^T(x) = h_\kappa^2(\sqrt{x_1}, \ldots,\sqrt{x_d})
   (1-|x|)^{\mu-1/2} /\sqrt{x_1 \cdots x_d},  
$$
where $\mu \ge 1/2$ and $h_\kappa$ is an reflection invariant weight function 
defined on $\RR^d$ and $h_\kappa$ is even in each of its variables. The case
$h_\kappa(x) = \prod_{i=1}^d|x_i|^{2\kappa_i}$ gives the classical weight
function on the simplex.

To emphasis the generality of the results, let us list several other families 
of weight functions beyond those in \eqref{eq:1.2}. For symmetric group of 
$d+1$ objects, 
\begin{equation}\label{eq:1.3}
h_\kappa(x) = \prod_{1 \le i,j \le d+1} |x_i - x_j|^{\kappa}, \qquad 
   \kappa \ge 0. 
\end{equation}
For hyperoctahedral group, the group generated by the reflections in $x_i=0$,
$1 \le i \le d+1$ and $x_i \pm x_j=0$, $1 \le i,j \le d+1$, 
\begin{equation}\label{eq:1.4}
h_\kappa(x) = \prod_{i=1}^{d+1} |x_i|^{\kappa_0} 
\prod_{1 \le i,j \le d+1} |x_i^2 - x_j^2|^{\kappa_1}, \qquad  \kappa_0, 
  \kappa_1 \ge 0. 
\end{equation}
Let us also mention that our results are often new even in the case of one 
dimension, $S^1$ or $[-1,1]$. For $S^1$ we have the dihedral group for which
\begin{equation}\label{eq:1.5}
h_\kappa(x_1,x_2) = |r^m \cos m \theta|^{\kappa_1}
  |r^m \sin m \theta|^{\kappa_2}, \qquad (x_1,x_2) = r (\cos \theta, 
    \sin \theta), 
\end{equation}
where $m$ is a positive integer, which can also be written in terms of 
Chebyshev polynomials of the first and the second kind. 

The paper is organized as follows. In the following section we give definitions
and discuss further properties and applications of the modulus of smoothness
defined in \cite{X02}. In Section 3 we define and discuss the modulus of 
smoothness of higher orders, use it to prove the direct and inverse theorems
and show that it is equivalent to a K-functional of higher order. The weighted
approximation on the unit ball $B^d$ and on the simplex $T^d$ is discussed in 
Section 4 and Section 5, respectively.

\section{means of $h$-harmonic expansion and modulus of smoothness} 
\setcounter{equation}{0}

\subsection{$h$-harmonic expansions}
Let $h_\kappa$ be the reflection invariant weight function defined in 
\eqref{eq:1.1}. The essential ingredient of the theory of $h$-harmonics is 
a family of first-order differential-difference operators, $\CD_i$, called 
Dunkl's operators, which generates a commutative algebra; these operators are 
defined by (\cite{D1})
$$
  \CD_i f(x) = \partial_i f(x) + \sum_{v \in R_+} k_v 
    \frac{f(x) -  f(x \sigma_v)} {\langle x, v\rangle}
        \langle v,\varepsilon_i\rangle, \qquad 1 \le i \le d+1, 
$$
where $\varepsilon_1, \ldots, \varepsilon_d$ are the standard unit vectors of 
$\RR^{d+1}$. The $h$-Laplacian is defined by $\Delta_h=\CD_1^2 + \ldots + 
\CD_{d+1}^2$ and it plays the role similar to that of the ordinary Laplacian. 
Let $\CP_n^{d+1}$ denote the subspace of homogeneous polynomials of degree $n$ 
in $d+1$ variables. An $h$-harmonic polynomial $P$ of degree $n$ is a 
homogeneous polynomial $P \in \CP_n^{d+1}$ such that $\Delta_h P  =0$. 
Furthermore, let $\CH_n^{d+1}(h_\kappa^2)$ denote the space of $h$-harmonic 
polynomials of degree $n$ in $d+1$ variables and define 
$$
 \langle f, g\rangle_\kappa : = a_\kappa \int_{S^d} f(x) g(x) 
    h^2_\kappa(x) d\omega(x). 
$$ 
Then $\langle P,Q\rangle_{\kappa} = 0$ for $P \in \CH_n^{d+1}(h_\kappa^2)$ and 
$Q \in \Pi_{n-1}^{d+1}$. The spherical $h$-harmonics are the restriction of 
$h$-harmonics on the unit sphere. Throughout this paper, we fix the value of 
$\lambda$ as 
\begin{equation}\label{eq:2.0}
  \lambda := \gamma_\kappa+ \frac{d-1}{2} \qquad \hbox{with} \qquad
   \gamma_\kappa =  \sum_{v\in R_+} \kappa_v.
\end{equation}
In terms of the polar coordinates $y =ry'$, $r = \|y\|$, the $h$-Laplacian
operator $\Delta_h$ takes the form (\cite{X01b})
$$ 
\Delta_h = \frac{\partial^2}{\partial r^2} + \frac{2 \lambda+1}{r} 
    \frac{\partial}{\partial r} + \frac{1}{r^2} \Delta_{h,0}, 
$$
where $\Delta_{h,0}$ is the (Laplace-Beltrami) operator on the sphere. Hence,
applying $\Delta_h$ to $h$-harmonics $Y \in \CH_n(h_\kappa^2)$ with $Y(y) = 
r^nY(y')$ shows that spherical $h$-harmonics are eigenfunctions of 
$\Delta_{h,0}$; that is,
\begin{equation}\label{eq:2.1}
 \Delta_{h,0} Y(x) = -n (n+ 2 \lambda) Y(x), 
      \qquad x \in S^d, \quad Y \in \CH_n^d(h_\kappa^2). 
\end{equation}
It is known that $\dim \CH_n(h_\kappa^2) =\dim \CP_n^d -
\dim \CP_{n-2}^d$ with $\dim \CP_n^d = \binom{n+d -1}{d}$. 

The standard Hilbert space theory shows that $L^2(h_\kappa^2) = 
\sum_{n=0}^\infty\bigoplus \CH_n^{d+1}(h_\kappa^2)$. That is, with each 
$f\in L^2(h_\kappa^2)$ we can associate its $h$-harmonic expansion
$$
  f(x) = \sum_{n=0}^\infty Y_n(h_\kappa^2;f,x), \qquad x \in S^d,  
$$ 
in $L^2(h_\kappa^2)$ norm. For the surface measure ($\kappa =0$), such a 
series is called the Laplace series (cf. \cite[Chapt. 12]{Er}). The orthogonal
projection $Y_n(h_\kappa^2): L^2(h_\kappa^2) \mapsto \CH_n^{d+1}(h_\kappa^2)$ 
takes the form 
\begin{equation}\label{eq:2.2}
 Y_n(h_\kappa^2;f,x):= 
   \int_{S^d} f(y) P_n(h_\kappa^2;x,y) h_\kappa^2(y) d\omega(y),
\end{equation} 
where $P_n(h_\kappa^2;x,y)$ is the reproducing kernel of the space of 
$h$-harmonics $\CH_n^{d+1}(h_\kappa^2)$. This kernel has a compact formula in
terms of the intertwining operator between the commutative algebra generated
by the partial derivatives and the one generated by Dunkl's operators. The 
intertwining operator $V_\kappa$ is a linear operator determined uniquely by
$$
V_\kappa \CP_n \subset \CP_n, \quad V_\kappa1=1, \quad 
\CD_i V_\kappa = V_\kappa \partial_i, \quad 1\le i\le d+1.
$$
 The compact formula of the reproducing kernel 
for $\CH_n^{d+1}(h_\kappa^2)$ is given by (\cite{X97b})
\begin{equation}\label{eq:2.3}
P_n(h_\kappa^2;x,y) = \frac{n+ \lambda}{\lambda} 
   V_\kappa [C_n^{\lambda} (\langle \cdot, y \rangle )](x), \qquad 
\end{equation}
where $C_n^\lambda$ is the usual Gegenbauer polynomial of degree $n$. 
If all $\kappa_v = 0$, $V_\kappa$ becomes the identity operator and the right 
hand is the so-called zonal harmonic. However, an explicit formula of 
$V_\kappa$ is known only in the case of symmetric group $S_3$ for three 
variables and in the case of the abelian group $\ZZ_2^{d+1}$. In the latter 
case, $V_\kappa$ is an integral operator given by (\cite{D2,X97a}) 
\begin{equation} \label{eq:2.4}
  V_\kappa f(x) = c_\kappa 
       \int_{[-1,1]^{d+1}} f(x_1 t_1, \ldots,x_{d+1} t_{d+1})
    \prod_{i=1}^{d+1} (1+t_i) (1-t_i^2)^{\kappa_i -1} d t,
\end{equation}
where $c_\kappa$ denotes the constant $c_\kappa =b_{\kappa_1} \ldots 
b_{\kappa_{d+1}}$ and $b_r^{-1} =\int_{-1}^1 (1-t^2)^{r-1}dt$. If some 
$\kappa_i =0$, then the formula holds under the limit relation
$$ 
 \lim_{\lambda \to 0} b_\lambda \int_{-1}^1 f(t) (1-t)^{\lambda -1} dt
  = [f(1) + f(-1)] /2. 
$$
This leads to an explicit formula for the reproducing kernel. One important 
property of the intertwining operator is that it is positive (\cite{Ros}); 
that is, $V p \ge 0$ if $p \ge 0$. Another important property is that 
$V_\kappa$ satisfies (\cite{X97b})
\begin{equation} \label{eq:2.4a}
a_\kappa \int_{S^{d}} V_\kappa f(\langle x,y\rangle) h_\kappa^2(y)d\omega(y) 
 = c_\lambda \int_{-1}^1 f(t) (1-t^2)^{\lambda -1/2}dt, \quad x \in S^d,  
\end{equation}
where $c_\lambda^{-1} = \int_{-1}^1 
(1-t^2)^{\lambda-1/2} dt = \Gamma(\lambda+1/2)\sqrt{\pi} /\Gamma(\lambda+1)$.

The compact formula \eqref{eq:2.3} of the reproducing kernel and the equation 
\eqref{eq:2.4a} indicate that the summability of the $h$-harmonic expansion is
related to Gegenbauer expansions. The Gegenbauer polynomials $C_n^\lambda(t)$ 
are orthogonal with respect to the weight function 
$$
  w_\lambda(t) = (1-t^2)^{\lambda -1/2}, \qquad -1 < t < 1.
$$ 
We denote by $\|g\|_{w_\lambda,p}$ the weighted $L^p(w_\lambda)$ norm for 
functions defined on $[-1,1]$,
$$
\|g\|_{w_\lambda,p} = \Big( c_\lambda \int_{-1}^1 |g(t)|^p w_\lambda(t) 
      dt\Big)^{1/p}
$$
for $1 \le p < \infty$ and $\|g\|_{w_\lambda,\infty} = \|g\|_\infty$ is the 
usual uniform norm on $[-1,1]$. For a function $g \in L^p(w_\lambda)$, its 
Gegenbauer expansion takes the form
$$
g(t) \sim \sum_{n=0}^\infty b_n \frac{n+\lambda}{\lambda} C_n^\lambda(t) \quad
\hbox{with} \quad b_n = \frac{c_\lambda}{C_n^\lambda(1)} \int_{-1}^1 g(t) 
     C_n^\lambda (t)w_\lambda(t) dt,
$$
since $\|C_n^\lambda\|_{w_\lambda,2}^2 = C_n^\lambda(1) \lambda/(n+\lambda)$
(cf. \cite[p. 80]{Szego}).

We will often encounter integral operators for functions on $S^d$ whose kernels
take the form of $G(x,y)= V_\kappa [g(\langle x, \cdot \,\rangle)](y)$, where 
$g: [-1,1]\mapsto \RR$. Such an operation defines a sort of convolution of 
the functions $f$ on $S^d$ and $g$ on $[-1,1]$, which we denote by 
$f \star_\kappa g$. Formally, we define

\begin{defn}
For $f \in L^p(h_\kappa^2)$ and $g \in L^1(w_\lambda;[-1,1])$, 
$$
(f\star_\kappa g)(x): = a_\kappa \int_{S^d} f(y) 
  V_\kappa[g(\langle x,\cdot\,\rangle)](y) h_\kappa^2(y) d\omega.
$$
\end{defn}

For the surface measure ($V_\kappa = id$), this is called spherical convolution
in \cite{CZ}. It satisfies many properties of the usual convolution in 
$\RR^{d+1}$. In particular, the familiar Young's inequality holds, which we 
state as follows.

\begin{prop}\label{prop:2.1}
Let $p,q,r \ge 1$ and $p^{-1} = r^{-1}+q^{-1}-1$. For $f \in L^q(h_\kappa^2)$ 
and $g \in L^r(w_\lambda;[-1,1])$, 
$$
\|f\star_\kappa g\|_{\kappa,p} \le \|f\|_{\kappa,q} \|g\|_{w_\lambda,r}. 
$$
\end{prop} 

\begin{proof}
The usual proof for Young's inequality works. We only need to notice that 
the inequality 
$$
\|G(x,\cdot)\|_{\kappa,r} \le \|g\|_{w_\lambda,r}, \qquad 
\hbox{where} \quad  G(x,y)= V_\kappa [g(\langle x, \cdot \,\rangle)](y), 
$$
holds, which we prove as follows: The fact that $V_\kappa$ is positive shows 
that $|Vg|\le V[|g|]$. Hence, the equation \eqref{eq:2.4a} leads to
$$
a_\kappa \int_{S^d} |G(x,y)| h_\kappa^2(y)d\omega \le 
a_\kappa \int_{S^d} V_\kappa [|g(\langle x, \cdot \,\rangle)|](y) 
  h_\kappa^2(y)d\omega 
= c_\lambda \int_{-1}^1 |g(t)| w_\lambda(t)dt.
$$
So that $\|G(x,\cdot)\|_{\kappa,1} \le  \|g\|_{w_\lambda,1}$. Evidently, we
also have $\|G(x,\cdot)\|_{\kappa,\infty} \le  \|g\|_{w_\lambda,\infty}$. 
The Riesz interpolation theorem shows then $\|G(x,\cdot)\|_{\kappa,r} \le 
\|g\|_{w_\lambda,r}$. 
\end{proof}

We will need the Ces\`aro $(C,\delta)$ means of the orthogonal expansion.
For $\delta >0$, the Ces\`aro $(C, \delta)$ means, 
$s_n^\delta$, of a sequence $\{c_n\}$ are defined by 
$$
s^\delta_n = (A_n^\delta)^{-1}\sum_{k=0}^n A_{n-k}^\delta c_k, 
  \qquad A_{n-k}^\delta = \binom{n-k+\delta}{n-k} 
$$
For the Gegenbauer expansion with respect to $w_\lambda$, we will use the 
notation
$$
P_n(w_\lambda;x,t) = \frac{n+\lambda}{\lambda} 
 \frac{C_n^\lambda(x)C_n^\lambda(t)}{C_n^\lambda(1)} \quad \hbox{and} \quad
P_n(w_\lambda;t) = P_n(w_\lambda;1,t). 
$$
Then the $(C,\delta)$ means of the Gegenbauer expansion, denoted by
$S_n^\delta(w_\lambda;f)$, can be written as an integral operator,
$$
S_n^\delta(w_\lambda;f,x) = c_\lambda \int_{-1}^1 f(t) 
 P_n^\delta(w_\lambda;x,t) w_\lambda(t) dt, 
$$
where, $P_n^\delta(w_\lambda;x,t)$ is the $(C,\delta)$ means of the sequence 
$\{P_n(w_\lambda;x,t)\}$. It is known that $S_n^\delta(w_\lambda;f)$ converges
to $f$ in $L^p(w_\lambda)$ norm if $\delta > \lambda$.

Using the notation of spherical convolution and \eqref{eq:2.3}, we can write 
\eqref{eq:2.2} as $Y_n(h_\kappa^2;f) = f\star_\kappa P_n(w_\lambda)$. 
Furthermore, we denote by $S_n^\delta(h_\kappa^2;f)$ 
the $(C,\delta)$ means of the sequence $\{Y_n(h_\kappa^2;f)\}$. Then, by 
\eqref{eq:2.3}, 
$$
 S_n^\delta(h_\kappa^2;f) = a_\kappa \int_{S^d} f(y) 
   V_\kappa[P_n^\delta(w_\lambda;\langle x,\cdot \,\rangle,1)](y) 
    h_\kappa^2(y)d\omega = f\star_\kappa P_n^\delta(w_\lambda),
$$
where we define $P_n^\delta(w_\lambda;t) = P_n^\delta(w_\lambda;t,1)$ and 
$\lambda$ as in \eqref{eq:2.0}. It is proved in \cite{X97b} that 
$S_n^\delta(h_\kappa^2;f)$ converges to $f$ in $L^p(h_\kappa^2)$ norm if 
$\delta > \lambda$, and the $S_n^\delta(h_\kappa^2; f)$ defines a positive 
operator if $\delta \ge 2\lambda +1$.
 
\subsection{Spherical means}
For the Lebesgue measure, the modulus of smoothness is defined via the 
spherical means, denoted by $T_\theta f$, 
\begin{equation} \label{eq:2.5}
 T_\theta f(x) = \frac{1}{\sigma_{d-1} (\sin \theta)^{d-1}}
     \int_{\langle x, y\rangle = \cos\theta} f(y) d\omega(y), 
\end{equation} 
where $\sigma_{d-1} = \int_{S^{d-1}} d\omega = 2 \pi^{d/2}/\Gamma(d/2)$, 
which is the surface area of $S^{d-1}$. The properties of the spherical means 
are well-known; see \cite{BBP,P}, for example. 

An extension of the spherical means associated with $h_\kappa^2 d\omega$, 
denoted by $T_\theta^\kappa f$, is defined in \cite{X02} in an indirect way:

\begin{defn}\label{defn1}
For $0 \le \theta \le \pi$, the means $T_\theta^\kappa$ are defined by 
\begin{equation}\label{eq:defn1}
c_\lambda \int_0^\pi T_\theta^\kappa f(x) g(\cos \theta) 
(\sin \theta)^{2\lambda} d\theta
= a_\kappa \int_{S^d} f(y) V_\kappa [g(\langle x,\cdot \rangle)](y) 
   h_\kappa^2(y)d\omega(y),
\end{equation} 
where $g$ is any $L^1(w_\lambda)$ function.
\end{defn} 

By Young's inequality in Proposition \ref{prop:2.1}, the right hand side
of \eqref{eq:defn1} is finite for $f \in L^1(h_\kappa^2)$ and $g \in 
L^1(w_\lambda)$ independent of $x$. Hence, for each $f \in L^1$ and each 
$x \in S^d$, the left hand side is a bounded linear functional on 
$L^1(w_\lambda)$ so that $T_\theta^\kappa f(x)$ is a unique function of 
$\theta$ in $L^\infty(w_\lambda)$ by \cite[Theorem 6.16, p. 127]{Rudin}. 

The definition is shown to be a proper 
extension of $T_\theta f$ in \cite{X02}, since the definition of $T_\theta f$ 
in \eqref{eq:2.5} satisfies \eqref{eq:defn1} when $\kappa =0$ and $V_\kappa = 
id$. Furthermore, if $f(x) = 1$ and $T_\theta^\kappa f(x) =1$, then 
\eqref{eq:defn1} is the equation \eqref{eq:2.4a}. By the definition of 
$f \star_\kappa g$, we can write \eqref{eq:defn1} as 
\begin{equation}\label{eq:defn2}
 (f \star_\kappa g)(x) = 
     c_\lambda \int_0^\pi T_\theta^\kappa f(x) g(\cos \theta) 
        (\sin \theta)^{2\lambda} d\theta.
\end{equation} 
The main properties of the means $T_\theta^\kappa f$ are given in the 
following proposition. 

\begin{prop}\label{prop:2.2}
The means $T_\theta^\kappa f$ satisfy the following properties:
\begin{enumerate}
\item Let $f_0(x) =1$, then $T_\theta^\kappa f_0(x) =1$. 
\item For $f\in L^1(h_\kappa^2)$, 
$$
  Y_n(h_\kappa^2; T_\theta^\kappa f) 
      = \frac{C_n^\lambda(\cos \theta)}{C_n^\lambda(1)} 
       Y_n(h_\kappa^2; f); 
$$
in particular, $\Delta_{h,0} T_\theta^\kappa f = T_\theta^\kappa(\Delta_{h,0} 
f)$ if $\Delta_{h,0}f \in L^1(h_\kappa^2)$. 
\item $T_\theta^\kappa: \Pi_n^d \vert_{S^d}\mapsto \Pi_n^d \vert_{S^d}$
and
$$
  T_\theta^\kappa f \sim \sum_{n=0}^\infty 
   \frac{C_n^\lambda(\cos \theta)}{C_n^\lambda(1)} Y_n(h_\kappa^2;f).  
$$
\item For $0 \le \theta \le \pi$,
$$
 T_\theta^\kappa f - f = \int_0^\theta (\sin s)^{-2\lambda} d s 
     \int_0^s T_t^\kappa(\Delta_{h,0} f) (\sin t)^{2\lambda} dt. 
$$
\item For $f\in L^p(h_\kappa^2)$, $1 \le p < \infty$, or $f\in C(S^d)$,
$$
 \|T_\theta^\kappa f\|_{\kappa,p} \le \|f\|_{\kappa,p} \qquad \hbox{and}
  \qquad \lim_{\theta \to 0} \|T_\theta^\kappa f -f \|_{\kappa,p} =0.
$$
\end{enumerate}
\end{prop}

\begin{proof}
These properties are analogous of those for the means $T_\theta$ (see 
\cite{BBP,P}). The first four properties were proved in \cite{X02}, and 
the property (5) was proved there only for $h_\kappa$ in \eqref{eq:1.2}. We
prove (5) for $h_\kappa^2$ in general below. By Young's inequality, 
$$
\left(a_\kappa \int_{S^d} \left| c_\lambda
  \int_0^\pi T_\phi^\kappa f(x) g(\cos \phi)
   (\sin \phi)^{2\lambda} d\phi \right|^p h_\kappa^2(x) d\omega \right)^{1/p}
  \le  \|f\|_{\kappa,p} \|g\|_{w_\lambda,1} 
$$
for any $g \in L^1(w_\lambda)$. For a fixed $\theta \in (0,\pi)$ and a 
sufficiently large integer $n$, choose $g$ such that $g(\cos \phi) = 1/B(n)$, 
$B_n(\theta) = c_\lambda \int_{\theta -1/n}^{\theta+1/n} (\sin \phi)^{2\lambda}
 d\phi$ if $|\phi - \theta|  \le 1/n$ and $g(\cos \theta) =0$ otherwise, it 
follows that $\|g\|_{w_\lambda,1} = 1$ and 
$$
\left(a_\kappa \int_{S^d} \left| \frac{1}{B_n(\theta)} 
 \int_{\theta-1/n}^{\theta+1/n}
   T_\phi^\kappa f(x) (\sin \theta)^{2\lambda} d\phi \right|^p h_\kappa^2(x) 
   d\omega \right)^{1/p}  \le  \|f\|_{\kappa,p}. 
$$
The fact that $T_\theta^\kappa f \in  L^\infty(w_\lambda)$ implies that 
$T_\theta^\kappa f \in L^1(w_\lambda)$. Hence, the inside integral converges to
$T_\theta^\kappa f(x)$ for almost all $\theta$ and the Fatou's lemma shows
$$
\|T_\theta^\kappa f\|_{\kappa,p}^p \le \liminf_{n \to \infty}
a_\kappa \int_{S^d} \left|\frac{1}{B_n(\theta)} 
 \int_{\theta-1/n}^{\theta+1/n}
   T_\phi^\kappa f(x) (\sin \theta)^{2\lambda} d\phi \right|^p h_\kappa^2(x) 
d\omega \le  \|f\|_{\kappa,p}^p
$$
for almost all $\theta$. Furthermore, the property (3) shows that for any 
polynomial $f$, $T_\theta^\kappa f$ is a continuous function on $\theta$. 
Hence, $\|T_\theta^\kappa f\|_{\kappa,p} \le c \|f\|_{\kappa,p}$ for all 
$\theta$ if $f$ is a polynomial. Taking, for example, $f_n$ as Ces\`aro 
$(C,\delta)$ means of the partial sum of $f$ with $\delta > 2 \lambda+1$, so 
that $f_n$ are positive and $f_n$ converge to $f$ in $L^p(h_\kappa^2)$ 
(\cite{X97b}). Then $Tf_n$ converges to $Tf$, so that $\|T_\theta^\kappa 
f\|_{\kappa,p} = \lim_{n\to \infty}  \|T_\theta^\kappa f_n\|_{\kappa,p} \le c 
\lim_{n\to\infty} \|f_n\|_{\kappa,p} = c\|f\|_{\kappa,p}$. This proves the 
first part of (5). The second part is a simple consequence of the first part
as in \cite{X02}.
\end{proof}
 

In the case $p =2$, (5) also follows from the Parseval identity. Indeed,
the property (2) implies that
$$
 \|T_\theta^\kappa f\|_{\kappa,2}^2 = \sum_{n=0}^\infty 
   \left|\frac{C_n^\lambda(\cos \theta)}{C_n^\lambda(1)}\right|^2  
      \|Y_n(h_\kappa^2; f)\|_{\kappa,2}^2 \le \sum_{n=0}^\infty 
        \|Y_n(h_\kappa^2; f)\|_{\kappa,2}^2 =  \|f\|_{\kappa,2}^2, 
$$
since $|C_n^\lambda(\cos \theta)| \le C_n^\lambda(1)$. For $h_\kappa^2$ in 
\eqref{eq:1.2}, the explicit formula of the intertwining operator $V_\kappa$ 
in \eqref{eq:2.4} is used to prove the inequality. 

One may attempt to define $\chi_\theta(\langle x,y\rangle)$ as the 
characteristic function of the set $S^d_\theta :=\{x \in S^{d}: 
\langle x,y\rangle = \cos \theta \}$, so that the ordinary spherical means 
$T_\theta f(x)$ can be written as 
$$
 T_\theta f(x) = \frac{1}{\sigma_{d-1} (\sin\theta)^{d-1}} 
   \int_{S^d} f(y) \chi_\theta(\langle x, y\rangle) d\omega(y).
$$
Since the set $S^d_\theta$ is a copy of $S^{d-1}$, which is one dimension 
lower than $S^d$, the function $\chi_\theta(t)$ behaves like a distribution 
function (Dirac delta function). Formally we can write 
$$
T_\theta^\kappa f(x)  = \frac{1} {c_\lambda (\sin\theta)^{2\lambda}}
a_\kappa \int_{S^d} f(y) V_\kappa[ \chi_\theta(\langle x, \cdot \rangle)]
  (y) h_\kappa^2(y)d\omega
$$
and use it to show property (5). However, since $\chi_\theta(t)$ is a 
distribution, we would have to show that $V \chi_\theta(\langle x, y \rangle)$ 
is a measure in order to carry through the calculation. Without further 
information on $V_\kappa$, this appears to be difficult. The author is in debt 
to an anonymous referee for pointing this out.

\subsection{Modulus of smoothness and best approximation}
The properties of $T_\theta^\kappa f$ are parallel to those of $T_\theta f$. 
They lead to the following definition of an analog of the modulus of 
smoothness:

\begin{defn} \label{defn2}
For $f \in L^p(h_\kappa^2)$, $1 \le p < \infty$, or  $f \in C(S^d)$, 
define
$$
\omega(f,t)_{\kappa,p} = \sup_{0<\theta \le t}
    \|T_\theta^\kappa f - f\|_{\kappa,p}.
$$
\end{defn}

Property (5) of Proposition \ref{prop:2.2} shows that, for $f \in 
L^p(h_\kappa^2)$,  $\omega(f;t)_{\kappa,p} \to 0$ if $t\to 0$. Moreover, 
it is proved in \cite{X02} that this modulus of smoothness is equivalent 
to a K-functional. Let 
$$
\CW^p(h_\kappa^2) =\{f \in L^p(h_\kappa^2): - k (k+2 \lambda)P_k(h_\kappa^2;f)
  = P_k(h_\kappa^2;g) \,\, \hbox{for some} \,\, g \in L^p(h_\kappa^2)\}.   
$$
Recall that $\Delta_{h,0}$ denote the spherical $h$-Laplacian. If 
$\Delta_{h,0}f \in L^p(h_\kappa^2)$, then $f \in \CW^p(h_\kappa^2)$ since we 
can take $g=\Delta_{h,0}f$ (see \eqref{eq:2.1}). The Peetre K-functional
between $L^p(h_\kappa^2)$ and $\CW^p(h_\kappa^2)$ is defined by 
$$
  K(f,t)_{\kappa,p} := \inf \big \{ \|f - g\|_{\kappa,p} + 
      t \|\Delta_{h,0}\, g\|_{\kappa,p}, \; g \in  \CW^p(h_\kappa^2) \big\}.
$$
The following proposition is proved in \cite{X02}:

\begin{prop} \label{prop2.4}
For $0 < t < \pi /2$, there exist constants $c_1$ and $c_2$ such that
$$
c_1 \omega(f,t)_{\kappa,p}  \le K(f,t^2)_{\kappa,p}
   \le c_2  \omega(f,t)_{\kappa,p}. 
$$ 
\end{prop}

Among the usual properties of the modulus of smoothness satisfied by 
$\omega(f,t)_{\kappa,p}$, we mention the following property: there is a
constant $c$ such that for $t > 0$, 
\begin{equation}\label{eq:2.7}
 \omega(f,  t \delta)_{\kappa,p} \le c \max \{t^2, 1\}
    \omega(f,\delta)_{\kappa,p}. 
\end{equation}
This important property is not obvious from the definition of 
$\omega(f,t)_{\kappa,p}$, it follows as a consequence of the equivalence 
to the K-functional. 

As an application of the definition of $\omega(f;t)_{\kappa,p}$ we consider the
convergence of the Poisson integral, defined by
$$
  P_r(f; x) = a_\kappa \int_{S^d} f(y) V_\kappa
  \Big[ \frac{1-r^2}{(1-2r \langle x, \cdot \, \rangle+ r^2)^{\lambda+1}}\Big]
   (y)h_\kappa^2(y)d\omega(y) 
$$ 
for $x \in S^d$ and $r < 1$. The kernel function of this integral is the 
Poisson kernel for the $h$-harmonics, which is equal to $\sum_{n=0}^\infty r^n 
P_n(h_\kappa^2,x,y)$ (\cite{D2}). 

\begin{prop} \label{prop:poisson}
For $f \in L^p(h_\kappa^2)$, $1 \le p < \infty$, or $f \in C(S^d)$, 
$p = \infty$,
$$
  \lim_{r \to 1-} \|P_r(f,x) - f(x)\|_{\kappa,p} = 0. 
$$
\end{prop}

\begin{proof}
The equation \eqref{eq:defn1} allows us to write
$$ 
P_r(f;x) = (f \star P_r)(x) \qquad \hbox{with} \quad P_r(t) = 
     \frac{1-r^2}{(1-2r \cos \theta + r^2)^{\lambda+1}}, 
$$
where $P_r(t)$ is equal to the Poisson kernel $P_r(w_\lambda;t,x)$ of the 
Gegenbauer expansion with $x =1$. In particular, $c_\lambda \int_0^\pi 
P_r(\cos\theta)(\sin \theta)^{2\lambda} d\theta =1$. Consequently, for any 
$\varepsilon > 0$, \eqref{eq:defn2} and the Minkowski inequality shows that 
\begin{align*}
 \|P_r f - f\|_{\kappa,p} & = 
 \Big \|c_\lambda \int_0^\pi (T_\theta^\kappa f(x)
    -f(x)) P_r(\cos \theta)
      (\sin\theta)^{2\lambda} d\theta \Big \|_{\kappa,p} \\
& \le c_\lambda \int_0^\pi \|T_\theta^\kappa f-f\|_{\kappa,p} 
    P_r(\cos \theta)(\sin\theta)^{2\lambda} d\theta\\
& \le \omega(f;\varepsilon)_{\kappa,p}+ 2 \|f\|_{\kappa,p} 
    c_\lambda \int_\varepsilon^\pi 
   P_r(\cos \theta) (\sin\theta)^{2\lambda} d\theta\\  
& \le \omega(f;\varepsilon)_{\kappa,p}+ 2 \|f\|_{\kappa,p} 
    \frac{1-r^2}{(1-2r \cos \varepsilon + r^2)^{\lambda+1}}. 
\end{align*}
Taking the limit $r\to 1$ gives $\|P_r f - f\|_{\kappa,p} \le 
\omega(f;\varepsilon)_{\kappa,p}$, which leads to the stated result since
$\varepsilon$ is arbitrary.
\end{proof}

In particular, if $f$ is continuous on $S^d$ then the proposition shows
that $P_r(f;x)$ converges uniformly to $f(x)$. One should compare this with 
\cite[Theorem 4.2]{X97a}, in which $P_r(f;x)$ is proved to converge 
pointwise to $f(x)$ for the case of $h_\kappa^2$ in \eqref{eq:1.2}; the proof 
there uses the explicit formula of $V_\kappa$ in \eqref{eq:2.4} and is much 
more involved. 

We can also use $\omega(f,t)_{\kappa,p}$ to study best approximation by 
polynomials and give the direct theorem and the inverse theorem. Although
more general theorems will be given in the following section, we give the
directed theorem below since its proof is constructive and further justifies 
the definition of $\omega(f,t)_{\kappa,p}$. 

\begin{thm}
For $f \in L^p(h_\kappa^2)$, $1 \le p \le \infty$,
$$
  E_n(f)_{\kappa,p} \le c\, \omega(f, n^{-1})_{\kappa,p}.
$$
\end{thm}

\begin{proof}
Let $s$ be a fixed positive integer and $s \ge \lambda + 2$. Let $n = m s$. 
We use the classical Jackson kernel defined by 
$$
   J_n(\cos \theta) = \frac{1}{\Omega_n} 
       \Big(\frac{\sin m\theta/2}{\sin \theta/2}\Big)^{2 s}, \qquad 
   0 \le \theta \le \pi,
$$
where $\Omega_n$ is a constant chosen such that 
$
c_\lambda \int_0^\pi J_n(\cos \theta) (\sin \theta)^{2 \lambda} d\theta = 1.
$
It is known that $J_n(t)$ is a nonnegative polynomial of degree $n = m s$
(cf. \cite[p. 203]{DL}). Using $J_n$ we consider $f\star_\kappa J_n$. 
Evidently, $f \star_\kappa J_n$ is a polynomial of degree at most $n$ and 
by \eqref{eq:defn1} we can write
$$
(f \star_\kappa J_n)(x) =  
  c_\lambda \int_0^\pi T_\theta^\kappa (f;x) J_n(\cos \theta)
    (\sin \theta)^{2 \lambda} d\theta.
$$
In particular, $(f_0 \star_\kappa J_n)(x) = 1$ for $f_0(x) =1$. Therefore, 
we have that 
\begin{align*}
\|f\star_\kappa J_n -f \|_{\kappa,p} & = c_\lambda 
  \Big \| \int_0^\pi J_n(\cos \theta) 
   (f(x) - T_\theta^\kappa (f;x)) (\sin \theta)^{2 \lambda} d\theta
  \Big \|_{\kappa,p} \\
& \le  c_\lambda \int_0^\pi \|f - T_\theta^\kappa (f)\|_{\kappa,p}
  J_n(\cos \theta)(\sin \theta)^{2 \lambda} d\theta
\end{align*}
Splitting the integral over $[0, \pi]$ into two integrals over $[0, 1/n]$ and
$[1/n,\pi]$, respectively, and using the definition of $\omega(f,t)$, we 
conclude that 
$$
\|f\star_\kappa J_n  - f\|_{\kappa,p} \le \omega(f,n^{-1})_{\kappa,p} +  
   c_\lambda \int_{1/n}^\pi \omega(f,\theta)_{\kappa,p}
    J_n(\cos \theta)(\sin \theta)^{2 \lambda} d\theta. 
$$
From the property \eqref{eq:2.7} of the modulus of continuity, it follows 
that for $\theta \ge n^{-1}$, 
$$
 \omega(f,\theta)_{\kappa,p} = \omega(f, n \theta/n )_{\kappa,p} \le 
  c \max\{1, n^2\theta^2\} \omega(f,n^{-1})_{\kappa,p} 
  = c n^2 \theta^2 \omega(f,n^{-1})_{\kappa,p}.
$$
Therefore, it follows that 
$$
\|f\star_\kappa J_n - f\|_{\kappa,p} \le \omega(f,n^{-1})_{\kappa,p} \left(1+
  c n^2 \int_{1/n}^\pi \theta^2 J_n(\cos \theta)(\sin \theta)^{2 \lambda} 
  d\theta \right).
$$
We write $A_n \asymp B_n$ if $c_1 \le |A_n/B_n| \le c_2$ for two constants 
$c_1$ and $c_2$ independent of $n$. As in \cite[p. 204]{DL}, we have 
\begin{align*}
\Omega_n & = \int_0^\pi \Big(\frac{\sin m\theta/2}{\sin \theta/2}\Big)^{2 s} 
   (\sin \theta)^{2 \lambda} d\theta
  \asymp \int_0^\pi \Big(\frac{\sin m\theta/2}{\theta}\Big)^{2 s} 
   \theta^{2 \lambda} (\cos \theta/2)^{2\lambda} d\theta \\
 & \asymp  m^{2s -2\lambda -1} \int_0^{m\pi/2} 
      \Big(\frac{\sin \phi}{\phi}\Big)^{2 s} 
   \phi^{2 \lambda} \Big(\cos \frac{\phi}{m}\Big)^{2\lambda} d\phi
  \asymp  m^{2s -2\lambda -1} 
\end{align*}
since $s \ge \lambda+2$. Similar estimate gives 
$$
\int_0^\pi \theta^2 J_n(\cos \theta) (\sin \theta)^{2\lambda} d\theta 
 = \frac{1}{\Omega_n} \int_0^\pi \theta^2 \Big(\frac{\sin m\theta/2}
     {\sin \theta/2}\Big)^{2 s} d\theta \asymp n^{-2}
$$
which proves the stated estimate. 
\end{proof}

We end this section with a proposition which will help us to deal with 
the spherical convolution $f \star_\kappa g$. 

\begin{prop} \label{prop:convolution} 
Assume $g \in L^1(w_\lambda,[-1,1])$ and $g(t) \sim \sum_{n=0}^\infty
b_n \frac{n+\lambda}{\lambda}C_n^\lambda(t)$. If $f\in L^p(h_\kappa^2)$
then the function $f \star_\kappa g$ is an element in $L^p(h_\kappa^2)$
and 
$$
  f \star_\kappa g \sim \sum_{n=0}^\infty b_n Y_n(h_\kappa^2;f).
$$
\end{prop}

\begin{proof}
The assumption on $g$ implies that $V_\kappa^{(y)} g(\langle x,y\rangle) \in 
L^1(h_\kappa^2)$ by \eqref{eq:2.4a}. Young's inequality with $r=1$ shows that 
$f \star_\kappa g \in L^p(h_\kappa^2)$. For any $Y_n^h\in \CH_n^d(h_\kappa^2)$,
part (2) of Proposition \ref{prop:2.2} gives
\begin{align*}
&  a_\kappa \int_{S^d} V_\kappa^{(y)}[g(\langle x,y\rangle)] Y_n^h(y)
    h_\kappa^2(y) d\omega = 
 \int_0^\pi T_\theta^\kappa Y_n^h (x) g(\cos \theta) (\sin \theta)^{2\lambda}
   d\theta\\
& \qquad  = \frac{c_\lambda}{C_n^\lambda(1)}
\int_{0}^\pi g(\cos \theta) C_n^\lambda(\cos\theta)(\sin \theta)^{2\lambda}
   d\theta  Y_n^h(x)  =b_n Y_n^h(x), 
\end{align*}
so that $V_\kappa^{(y)}g(\langle x, y \rangle) \sim \sum_{n=0}^\infty b_n 
P_n(h_\kappa^2;x,y)$. Let $P_r(f)$ denote the Poisson integral of $f$. For
$0 \le r <1$, the expansion
$$
 P_r(f,x) = \sum_{n=0}^\infty r^n Y_n(h_\kappa^2;f,x) 
$$
holds uniformly for $x \in S^d$, which allows us to integrate term by 
term. Consequently, for $0<r<1$ the function $P_rf \star_\kappa g$ is well 
defined and it follows that
$$\langle P_rf \star_\kappa g,Y_n^h \rangle_{\kappa,p} = 
  b_n r^n \langle f, Y_n^h\rangle_{\kappa,p}$$ 
for any $Y_n^h \in \CH_n^d(h_\kappa^2)$. Furthermore, Young's inequality 
shows that 
$$
\|P_rf \star_\kappa g - f\|_{\kappa,p} 
    \le \|P_r(f)- f\|_{\kappa,p} \|g\|_{w_\lambda,1}
$$
so that $\|P_rf \star_\kappa g  - f \star_\kappa g \|_{\kappa,p} \to 0$ as 
$r \to 1$ by Proposition \ref{prop:poisson}. Consequently, 
$\langle f \star_\kappa g, Y_n^h\rangle_{\kappa,p} = 
\lim_{r \to 1} \langle P_rf \star_\kappa g, Y_n^h\rangle_{\kappa,p} = 
b_n \langle f, Y_n^h\rangle_{\kappa,p}$.
\end{proof}

\section{Modulus of smoothness and K-functional of higher order} 
\setcounter{equation}{0}

If $s$ is an integer, then $(-\Delta_{h,0})^s = -\Delta_{h,0} 
(-\Delta_{h,0})^{s-1}$ is well defined. By \eqref{eq:2.1}, $(-\Delta_{h,0})^s 
g \sim \sum_{n=1}^\infty (n(n+2\lambda))^s Y_n(h_\kappa^2;g)$. We can also 
define the fractional power of $\Delta_{h,0}$ using the $h$-harmonic 
expansion. 

\begin{defn}
Let $r$ be a positive integer. Define $(-\Delta_{h,0})^{r/2}g$ by 
$$
(-\Delta_{h,0})^{r/2} g 
 \sim \sum_{n=1}^\infty (n(n+2\lambda))^{r/2} Y_n(h_\kappa^2;g).
$$   
Furthermore, define the function space $\CW_r^p(h_\kappa^2)$ by
$$
\CW^p_r(h_\kappa^2) =\{f \in L^p(h_\kappa^2): (k (k+2 \lambda))^{\frac{r}2}
 P_k(h_\kappa^2;f) = P_k(h_\kappa^2;g) \,\, \hbox{for some} \,\, g \in 
  L^p(h_\kappa^2)\}.
$$
\end{defn}

By definition, the space $\CW^p(h_\kappa^2)$ defined
in the previous section is the same as $\CW^p_2(h_\kappa^2)$. 
If $(-\Delta_{h,0})^{r/2}f \in L^p(h_\kappa^2)$ then $f \in 
\CW_r^p(h_\kappa^2)$ since we can take $g = (-\Delta_{h,0})^{r/2}f$ in the 
definition. On the other hand, if $f \in \CW_r^p(h_\kappa^2)$, then 
\eqref{eq:2.1} shows that $g$ and
$(-\Delta_{h,0})^{r/2}f$ have the same coefficients in their $h$-harmonic 
expansions, so that $g = (-\Delta_{h,0})^{r/2}f$ in the $L^p(h_\kappa^2)$ 
norm, which shows that $(-\Delta_{h,0})^{r/2}f \in L^p(h_\kappa^2)$.  
In fact, if $F$ is a function whose coefficients $\langle F,Y_n^h
\rangle_{\kappa,p} =0$ for all $Y_n^h \in \CH_n^{d+1}(h_\kappa^2)$, then 
considering the Ces\`aro $(C, \lambda+1)$ means (so that the means converge) 
shows that $\|F\|_{\kappa,p} =0$.

For each $r > 0$, by a result of \cite[Theorem 1 and Theorem 3]{AW}, there is
a function $\phi_r(x)$ such that $\phi_r$ is continuous on $[-1,1)$, $\phi_r
\in L^1(w_\lambda, [-1,1])$, and 
$$
   \phi_r(t) \sim \sum_{n=1}^\infty (n(n+2\lambda))^{-r/2} 
      \frac{n+\lambda}{\lambda} C_n^\lambda(t). 
$$
Furthermore, it follows from the theorems in \cite{AW} that $\phi_r \in
L^p(w_\lambda,[-1,1])$ if $r > (2\lambda+1)/q$, where $p^{-1}+q^{-1}=1$, and
$\phi_r$ is continuous on $[-1,1]$ if $r > 2\lambda+1$. 

\begin{lem} \label{lem:3.2}
If $f \in \CW_p^r (h_\kappa^2)$ then $f(x) = (-\Delta_{h,0})^{r/2} f 
\star_\kappa \phi_r$ in $L^p(h_\kappa^2)$. 
\end{lem}

\begin{proof}
For $f \in L^p(h_\kappa^2)$, the Proposition \ref{prop:convolution} shows that 
$f \star_\kappa \phi_r \in L^p (h_\kappa^2)$ and 
$$
 f \star_\kappa \phi_r 
   \sim \sum_{n=0}^\infty (n(n+2\lambda))^{-r/2} Y_n(h_\kappa^2;f).
$$  
If $f\in \CW_p^r(h_\kappa^2)$, then $(-\Delta_{h,0})^{r/2}f\in L^p(h_\kappa^2)$
so that $(-\Delta_{h,0})^{r/2} f \star_\kappa \phi_r \in L^p(h_\kappa^2)$ and 
$$
(-\Delta_{h,0})^{r/2} f \star_\kappa \phi_r
   \sim \sum_{n=1}^\infty  (n(n+2\lambda))^{-r/2} 
 Y_n(h_\kappa^2; (-\Delta_{h,0})^{r/2} f) = \sum_{n=1}^\infty 
 Y_n(h_\kappa^2;f) \sim f.
$$   
It follows that $(-\Delta_{h,0})^{r/2} f \star_\kappa \phi_r= f$ in 
$L^p(h_\kappa^2)$. 
\end{proof}

In the following we shall use $c$ to denote a generic constant, which depends 
only on $d$, $p$, $r$ and $\kappa$ and whose value may be different from line 
to line.

\begin{thm} \label{thm:3.2}
For $f \in \CW_p^r (h_\kappa^2)$, $1 \le p \le \infty$,
$$
 E_n(f)_{\kappa,p} \le c\, n^{-r} \|(-\Delta_{h,0})^{-r/2} f \|_{\kappa,p}.
$$
\end{thm}
 
\begin{proof}
Let $\sigma$ be a positive integer, $\sigma > 2 \lambda+1$ so that 
$P_n^\sigma (w_\lambda;x,t)$ is nonnegative. Using summation by parts 
repeatedly on the expansion of $\phi_r(t)$, we can write 
$$
 \phi_r(t) = \sum_{k=0}^\infty \Delta^{\sigma+1} \mu(k) \binom{k+\sigma}{k}
     P_k^\sigma(w_\lambda;t,1), \qquad \mu(k) = (k(k+2\lambda))^{-r/2}, 
  \quad k \ge 1 
$$
and $\mu(0) =0$, where $\Delta^m \mu(k)$ denotes the $m$-th order finite 
difference, defined by $\Delta \mu(t) = \mu(t) - \mu(t+1)$ and $\Delta^{m+1}
= \Delta (\Delta^m)$. Let $q_n$ be the $n$-th partial sum of the above series. 
Then $(-\Delta_{h,0})^{r/2} f \star_\kappa q_n$ is evidently a polynomial of 
degree at most $n$. It follows from the Lemma \ref{lem:3.2} that 
\begin{align*}
\|f(x) - (-\Delta_{h,0})^{r/2} f \star_\kappa q_n \|_{\kappa,p} 
 & = \|(-\Delta_{h,0})^{r/2} f \star_\kappa (\phi_r - q_n)\|_{\kappa,p} \\
 & \le  \|(-\Delta_{h,0})^{r/2} f\|_{\kappa,p}  
   \|\phi_r(\cos \theta) - q_n(\cos\theta)\|_{w_\lambda,1}.
\end{align*}
Since $P_n^\sigma(w_\lambda)$ is nonnegative, $\|P_n^\sigma(w_\lambda;\cdot,1)
\|_{\kappa,1} =1$. The finite difference satisfies $\Delta^m \mu(t) = (-1)^m
\mu^{(m)}(\xi)$ for some $\xi$ between $t$ and $t + m$. With 
$\mu(t) = (t(t+2\lambda))^{-r/2}$, it is easy to verify that 
$$
 \mu^{(\sigma+1)}(t) = (-1)^{\sigma+1}r(r+1) \cdots (r+\sigma) \mu(t)  
    t^{-\sigma-1},
$$ 
so that $|\Delta^{\sigma+1} \mu(k)| \le c k^{- r - \sigma -1}$. Therefore,
using $\binom{k+\sigma}{k} \sim k^{\sigma}$, it follows that 
\begin{align*}
 \|\phi_r(\cos \theta) - q_n(\cos\theta)\|_{w_\lambda,1}
 & = \Big\|\sum_{k=n+1}^\infty \Delta^{\sigma+1} \mu(k) \binom{k+\sigma}{k}
    P_k^\sigma(w_\lambda;t,1)\Big\|_{w_\lambda,1} \\
 & \le c \sum_{k=n+1}^\infty |\Delta^{\sigma+1} \mu(k)| \binom{k+\sigma}{k}\\
 &  \le c \sum_{k=n+1}^\infty k^{-r -1} \le c n^{-r},
\end{align*}
which completes the proof. 
\end{proof}

For the ordinary spherical harmonics the above theorem was essentially proved
in \cite{Kam}, and the proof follows the same line there. The idea of using 
summation by parts and the Ces\`aro means plays an important role in the 
proof of the direct theorem for the Lebesgue measure in \cite{Rus} which we 
follow in the development below.

Using the space of $\CW_p^r(h_\kappa^2)$ we can define a K-functional of 
$r$-th order. 

\begin{defn} 
For $r \in \NN$, the K-functional between $L^p(h_\kappa^2)$ and 
$\CW_r^p(h_\kappa^2)$ is 
$$
K_r(f;t)_{\kappa,p} := \inf \big \{ \|f - g\|_{\kappa,p} + t^r 
 \|(-\Delta_{h,0})^{r/2}\, g\|_{\kappa,p}, \; g \in \CW_r^p(h_\kappa^2) \big\}.
$$
\end{defn}

We note that the relation with the K-functional defined in the previous 
section is $K(f,t^2)_{\kappa,p} = K_2(f,t)_{\kappa,p}$. Since 
$E_n(f)_{\kappa,p}  \le \|f-g\|_{\kappa,p} + E_n(g)_{\kappa,p}$, it follows
from Theorem \ref{thm:3.2} that 
\begin{equation}\label{eq:3.0}
E_n(f)_{\kappa,p} \le c K_r(f;n^{-1})_{\kappa,p}.
\end{equation}
Our goal is to prove the similar result using the modulus of smoothness 
of higher order, which is defined using the power of the operator 
$I- T_\theta^\kappa$. If $s$ is an integer, then $(I-T_\theta^\kappa)^s = 
\sum_{k=0}^s (-1)^s \binom{s}{k}(T_\theta^\kappa)^k$. We can define the
fractional power of the operator similarly. However, recall that 
$T_\theta^\kappa Y_n(h_\kappa^2;f)= (C_k^\lambda(\cos\theta)/
C_k^\lambda(1))Y_n(h_\kappa^2;f)$ by the Proposition \ref{prop:2.2}, 
we define the fractional power using the $h$-harmonic expansion. 

\begin{defn}
Let $r \in \NN$. Define 
$$
(I- T_\theta^\kappa)^{r/2} f \sim \sum_{n=0}^\infty 
     (1- R_n^\lambda(\cos \theta))^{r/2} Y_n(h_\kappa^2;f), \qquad
 R_k^\lambda(t):= C_k^\lambda(t)/C_k^\lambda(1).
$$
For $f\in L^p(h_\kappa^2)$, $1 \le p < \infty$, or 
$f \in C(S^d)$, define
$$
\omega_r(f,t)_{\kappa,p} 
  := \sup_{0\le \theta \le t} \|(I-T_\theta^\kappa)^{r/2}\|_{\kappa,p}.
$$
\end{defn}

We note that the modulus of smoothness defined in the previous section 
corresponds to $\omega_2(f,t)_{\kappa,p}$. For the Lebesgue measure 
($\kappa =0$), such a definition was given in \cite{Rus} and the case $r$ 
being an even integer had appeared in several early references (see the 
discussion in \cite{Rus}). Some of the properties of 
$\omega_r(f,t)_{\kappa,p}$ is collected below.

\begin{prop}
The modulus of smoothness $\omega_r(f,t)_{\kappa,p}$ satisfies:
\begin{enumerate}
\item $\omega_r(f,t)_{\kappa,p} \to 0$ if $t \to 0$;
\item $\omega_r(f,t)_{\kappa,p}$ is monotone nondecreasing on $(0,\pi)$;
\item $\omega_r(f+g,t)_{\kappa,p} \le \omega_r(f,t)_{\kappa,p}+
         \omega_r(g,t)_{\kappa,p}$; 
\item For $0 < s < r$, 
$$
\omega_r(f,t)_{\kappa,p} \le 2^{[(r-s+1)/2]} \omega_s(f,t)_{\kappa,p};
$$
\item If $(-\Delta_{h,0})^k f \in L^p(h_\kappa^2)$, $k \in \NN$, then 
for $r > 2 k$
$$
\omega_r(f,t)_{\kappa,p} 
 \le c\, t^{2k} \omega_{r-2k}((-\Delta_{h,0})^k f,t)_{\kappa,p}.
$$
\end{enumerate}
\end{prop}

\begin{proof}
If $\|f\|_{\kappa,p} < 1$ then $\|T_\theta^\kappa f\|_{\kappa,p} < 1$ and
$$
 (I- T_\theta^\kappa)^{r/2} f = \sum_{n=0}^\infty
   (-1)^n \binom{r/2}{n}  (T_\theta^\kappa)^n f  
$$
in the $L^p(h_\kappa^2)$ norm, so that $\|(I- T_\theta^\kappa)^{r/2} f
\|_{\kappa,p} \le 2^{[(r+1)/2]} \|f\|_{\kappa,p}$. Notice also that the space 
of polynomials is dense in $L^p(h_\kappa^2)$. The first four properties of 
$\omega_r(f,t)_{\kappa,p}$ follow easily from the definition and the above 
inequality. To prove the fifth 
property, we notice that $(I-T_\theta^\kappa)^{r/2}$ commutes with 
$\Delta_{h,0}$ by definition. Property (4) of Proposition
\ref{prop:2.2} shows that $\|(I-T_\theta^\kappa)f\|_{\kappa,p} \le c t^2
\|(-\Delta_{h,0})f\|_{\kappa,p}$, using this inequality repeatedly gives the 
stated property.
\end{proof}

We will prove a directed theorem using $\omega_r(f,t)_{\kappa,p}$ and also 
prove that $\omega_r(f,t)_{\kappa,p}$ is equivalent to $K_r(f;t)_{\kappa,p}$.
For the Lebesgue measure ($\kappa =0$), this was studied by several authors
and finally solved by Rustamov in \cite{Rus} (also see \cite{Rus} for the 
historical account). In the remainder of this section, we will follow the 
approach in \cite{Rus} closely. Along the way, we will point out the 
similarity and the major difference in the proof. 

Let $\eta \in C^\infty[0,+\infty)$ be a function defined by $\eta(x) =1$ for
$0 \le x \le 1$ and $\eta(x) =0 $ if $x \ge 2$. Define a sequence of operators
$\eta_n$ for $n \in \NN$ by 
$$
 \eta_n f: = \sum_{k=0}^\infty \eta\Big(\frac{k}{n}\Big) Y_k(h_\kappa^2;f)
   =  f \star_\kappa \eta_n(w_\lambda), \qquad 
  \eta_n(w_\lambda, t) = \sum_{n=0}^\infty \eta\Big(\frac{k}{n}\Big) 
       P_k(w_\lambda;t).
$$
Since $\eta(k/n) = 0$ if $k \ge 2n$, the infinite series terminates at 
$k = 2n-1$ so that $\eta_n$ is a spherical polynomial of degree at most 
$2n-1$. Furthermore, if $P$ is a spherical polynomial of degree at most $n$, 
then $Y_k(h_\kappa^2;P)=0$ for $k > n$ and the definition of $\eta$ shows 
that $\eta_n P = P$. The main properties of $\eta_n$ are given in the 
following lemma.

\begin{prop} \label{prop:3.7}
Let $f \in L^p(h_\kappa^2)$, $1 \le p \le \infty$, then
\begin{enumerate}
\item $\eta_n f \in \Pi_{2n-1}^{d+1}$ and $\eta_n P  = P$ for 
$P\in \Pi_n^{d+1}$; 
\item for $n \in \NN$, $\|\eta_n f \|_{\kappa,p} \le c \|f\|_{\kappa,p}$; 
\item for $n \in \NN$, 
$$
 \|f- \eta_n f \|_{\kappa,p} \le c E_n (f)_{\kappa,p}.
$$
\end{enumerate}
\end{prop}

\begin{proof} 
The proof uses the summation by parts and $(C,\delta)$ means as in the 
proof of Theorem \ref{thm:3.2}. Young's inequality (with $r=1$) gives
$
\|\eta_n f \|_{\kappa,p} \le  \|f\|_{\kappa,p} 
  \|\eta_n(w_\lambda)\|_{w_\lambda,1}  
$
and, with the notation as in the proof of Theorem \ref{thm:3.2}, we can write 
\begin{align*}
\|\eta_n (w_\lambda)\|_{w_\lambda,1} 
 & = \Big \|\sum_{k=1}^\infty \Delta^{\sigma +1} \eta\Big(\frac{k}{n}\Big)  
      \binom{k+\sigma}{k} P_k^\sigma(w_\lambda;t,1) \Big\|_{w_\lambda,1} \\
 & \le c \sum_{k=1}^{2n} \Big|\Delta^{\sigma +1} \eta\Big(\frac{k}{n}\Big)
   \Big| k^\sigma   \le c
\end{align*}  
since $\eta \in C^\infty[0,+\infty)$ implies that $|\Delta^{\sigma +1} 
\eta(k/n)| \le c n^{-\sigma-1}$. This proves the part (2). The part
(3) is an easy consequence of (1), (2) and triangle inequality.
\end{proof}

For Lebesgue measure ($\kappa =0$), this construction appears in \cite{Rus}
and the proof is identical. We repeated the proof since it is simple and
reinforces the idea used in the proof of Theorem \ref{thm:3.2}. The same idea 
is also used in the following proposition, whose proof essentially follows 
from the analogous result for the Lebesgue measure in \cite{Rus}. We give an
outline of the proof.  

\begin{prop} \label{prop:3.8}
Suppose $0 < t < \pi/(2(n+\lambda))$. For $f\in L^p(h_\kappa^2)$,
$1 \le p \le \infty$, 
\begin{equation} \label{eq:3.1}
\|(-\Delta_{h,0})^{r/2} \eta_n f \|_{\kappa,p} \le c \,t^{-r}
   \|(I-T_t^\kappa)^{r/2} f \|_{\kappa,p}; 
\end{equation} 
furthermore, 
\begin{equation} \label{eq:3.2}
\|(I-T_t^\kappa)^{r/2}\eta_n f \|_{\kappa,p} \le c \,t^r
   \|(-\Delta_{h,0})^{r/2} f \|_{\kappa,p}. 
\end{equation} 
\end{prop}

\begin{proof}
Since $P_n(h_\kappa^2;f) = f \star_\kappa P_n(w_\lambda)$, it follows from
\eqref{eq:2.1} and Proposition \ref{prop:2.2} that 
\begin{align*}
 (-\Delta_{h,0})^{r/2} \eta_n f & = \sum_{k=1}^\infty \eta\Big(\frac{k}{n}\Big)
    (k(k+2\lambda))^{r/2} Y_n(h_\kappa^2;f) \\
  & = \sum_{k=1}^\infty \eta\Big(\frac{k}{n}\Big)
    \alpha_\theta(k) Y_n(h_\kappa^2; (I-T_\theta^\kappa)^{r/2}f) 
\end{align*}
where $\alpha_\theta(k) =  (k(k+2\lambda))^{r/2}/
 (1-R_k^\lambda(\cos \theta))^{r/2}$. Therefore, we have 
$$
(-\Delta_{h,0})^{r/2} \eta_n f = (I-T_\theta^\kappa)^{r/2}f \star_\kappa 
   \alpha(\theta), 
$$
where, using summation by parts as in the proof of Theorem \ref{thm:3.2},
$$
\alpha(\theta;t) := \sum_{k=1}^\infty 
  \eta\Big(\frac{k}{n}\Big)\alpha_\theta(k) P_k(w_\lambda;t)
  = \sum_{k=1}^\infty \Delta^{\sigma +1} \Big[\eta\Big(\frac{k}{n}\Big)
  \alpha_\theta(k) \Big]  
   \binom{k+\sigma}{k} P_k^\sigma(w_\lambda;t).
$$
The same consideration also shows that 
$$
(I-T_\theta^\kappa)^{r/2}\eta_n f = (-\Delta_{h,0})^{r/2}f \star_\kappa 
   \beta(\theta), 
$$
where, defining $\beta_\theta(k) = (1-R_k^\lambda(\cos \theta))^{r/2}/
(k(k+2\lambda))^{r/2}$, we have
$$
\beta(\theta;t) := \sum_{k=1}^\infty 
  \eta\Big(\frac{k}{n}\Big)\beta_\theta(k) P_k(w_\lambda;t)
  = \sum_{k=1}^\infty \Delta^{\sigma +1} \Big[\eta\Big(\frac{k}{n}\Big)
  \beta_\theta(k) \Big]  
   \binom{k+\sigma}{k} P_k^\sigma(w_\lambda;t).
$$
Thus, as in the proof of Proposition \ref{prop:3.7}, the proof of the 
two inequalities is reduced to prove 
\begin{align*}
 \sum_{k=1}^{2n-1} \Big|\Delta^{\sigma +1} \Big[\eta\Big(\frac{k}{n}\Big)
     \alpha_\theta(k) \Big]\Big| k^\sigma   \le c \quad \hbox{and} \quad
 \sum_{k=1}^{2n-1} \Big|\Delta^{\sigma +1} \Big[\eta\Big(\frac{k}{n}\Big)
     \beta_\theta(k) \Big]\Big| k^\sigma   \le c.
\end{align*}
These inequalities are proved in \cite{Rus}. Using the Leibniz rule for 
finite difference and the relation between finite differences and derivatives, 
the main task is to estimate the derivatives of $\alpha_\theta$ and 
$\beta_\theta$, where the formula 
$$
 R_k(\cos\theta) = c_\lambda(\sin \theta)^{1-2\lambda} 
   \int_0^1 [\cos \lambda \phi - \cos (k+\lambda)\phi]
      (\cos \phi - \cos \theta)^{\lambda-1} d\phi 
$$
is used to extend the definition of $R_k(\cos \theta)$ as a function of 
$k \in \RR$. The estimates are rather involved, see \cite{Rus} for details. 
\end{proof}

Since $\eta_n f$ preserves polynomials of degree $n$, one immediate 
consequence of the inequality \eqref{eq:3.1} is the following inequality:

\begin{prop} \label{prop:3.9}
Suppose $0 < t < \pi/(2(n+\lambda))$. For any polynomial $P_n \in \Pi_n^{d+1}$,
\begin{equation*} 
\|(-\Delta_{h,0})^{r/2} P_n \|_{\kappa,p} \le c \,t^{- r}
   \|(I-T_t^\kappa)^{r/2} P_n \|_{\kappa,p}. 
\end{equation*} 
\end{prop}

The inequality in the proposition is called an inequality of
Riesz-Bernstein-Nikolskii-Stechkin type in \cite{Rus}. It implies, in
particular, the Bernstein type inequality 
\begin{equation}\label{eq:3.3} 
\|(-\Delta_{h,0})^{r/2} P_n \|_{\kappa,p} \le c \, n^{r}\|P_n \|_{\kappa,p}. 
\end{equation} 

Another important consequence of the inequality \eqref{eq:3.1} is the 
following:

\begin{prop} \label{prop:3.10}
For $f \in L^p(h_\kappa^2)$, $1 \le p \le \infty$, 
\begin{equation*} 
 \|(-\Delta_{h,0})^{r/2} \eta_n f\|_{\kappa,p} \le c \,n^r
   \omega_r(f, \pi /(2(n+\lambda)))_{\kappa,p}. 
\end{equation*} 
\end{prop}

This corollary will help us to prove a direct theorem (Jackson type estimate).
Much of the difficulty of the proof comes from the fact that the relation
$$
 \omega_r(f;\delta t)_{\kappa,p} \le  c\, \max\{1,\delta^r\} 
     \omega_r(f;t)_{\kappa,p}
$$
is not obvious from the definition. It will be established as a corollary of
the equivalence between $\omega_r(f;t)_{\kappa,p}$ and $K_r(f;t)_{\kappa,p}$
(see Corollary \ref{cor:3.14}). However, the proof of one direction of the 
equivalence will use the Jackson type theorem, where the other direction is 
needed in the proof of the Jackson type theorem. This explains why only half 
of the equivalence is given in the following proposition.

\begin{prop} \label{prop:3.11}
For $f \in L^p(h_\kappa^2)$, $1 \le p \le \infty$, 
$\omega_r(f;t)_{\kappa,p} \le c K_r(f;t)_{\kappa,p}$.
\end{prop}

\begin{proof}
Let $g \in \CW_r^p(h_\kappa^2)$. Using \eqref{eq:3.2} we get
\begin{align*}
\omega_r(f;t)_{\kappa,p} & \le c \|f-g\|_{\kappa,p} + 
   \omega_r(g;t)_{\kappa,p} \\    
 & \le c ( \|f-g\|_{\kappa,p} + t^r \|(-\Delta_{h,0})^{r/2} g\|_{\kappa,p}.
\end{align*}
Taking infimum over $g \in \CW_r^p(h_\kappa^2)$ gives the stated inequality.
\end{proof}

We are ready to prove the Jackson type estimate. The proof follows the 
one given in \cite{Rus} for the Lebesgue measure, but it differs in a major
step: the proof in \cite{Rus} uses a lemma that is established by a complicated
limit argument (\cite[Lemma 3.9]{Rus}), which is in fact being questioned
in \cite{LW}. Our proof does not depend on such a lemma.

\begin{prop} \label{prop:3.12}
For $f \in L^p(h_\kappa^2)$, $1 \le p \le \infty$, 
$$
   E_n(f)_{\kappa,p} \le c\, \omega_r(f;\pi/(2(n+\lambda)))_{\kappa,p}.
$$
\end{prop}

\begin{proof}
We can assume that $f$ is orthogonal to constants with respect to $h_\kappa^2 
d\omega$ since the constant term has no impact on the best approximation or 
the value of $\omega_r(f;t)_{\kappa,p}$. We define a sequence $n_j$, 
$j=0,1,2,\ldots$, as follows:
$$
  n_0 =1,\qquad n_{j+1} = \inf\{n: \omega_r(f;\pi/(2(2n+\lambda))) \le 
       \omega_r(f;\pi/(2(2n_j+\lambda)))_{\kappa,p}/2 \}
$$
for $j \ge 0$. The fact that $\omega_r(f;t)_{\kappa,p}$ is monotone 
nondecreasing on $(0,\pi)$ shows that $n_j \to \infty$ as $j \to \infty$. 
We claim that the sequence $n_k /n_{k-1}$ is bounded by a constant $c$.
Suppose otherwise; then there is a subsequence that goes to infinity. Without
lose of generality, we can assume that $n_k /n_{k-1}$ goes to infinity as
$k \to \infty$. By the definition of $n_k$, 
$$
   \omega_r(f;\pi/(2(2 n_k -2+\lambda))) > 
       \omega_r(f;\pi/(2(2n_{k-1}+\lambda)))_{\kappa,p}/2.
$$
Write $\delta_k = \pi/(2(2n_k-2+\lambda))$. Then using the above inequality,
Proposition \ref{prop:3.10} and Proposition \ref{prop:3.11} we have
\begin{align*}
 n_{k-1}^{-r} \|(-\Delta_{h,0})^{r/2} \eta_{2n_{k-1}} f\|_{\kappa,p}
   &\le c \,\omega_r(f; \pi/(2(2 n_{k-1} +\lambda))) \\
   & \le c \,\omega_r(f; \delta_k)_{\kappa,p}  
    \le c \, K_r(f;\delta_k)_{\kappa,p}. 
\end{align*}
Consequently, for any given $\varepsilon >0$, the definition 
of $K_r(f;t)_{\kappa,p}$ shows that for some $g\in \CW_r^p(h_\kappa^2)$,  
$$
 \|(-\Delta_{h,0})^{r/2} \eta_{2n_{k-1}} f\|_{\kappa,p}
  \le c\, n_{k-1}^{r} K_r(f; \varepsilon)_{\kappa,p}      
  + c (n_{k-1}/n_k)^r \|(-\Delta_{h,0})^{r/2}g \|_{\kappa,p}.      
$$
Let $\varepsilon$ go to zero and then let $k \to \infty$. Since 
$n_{k-1}/n_k \to 0$, Fatou's lemma shows that 
$\|(-\Delta_{h,0})^{r/2} f\|_{\kappa,p} =0$. This shows, however, that 
$f$ is a constant, which is a contradiction to our assumption. 

By Proposition \ref{prop:3.7}, $f = \sum_{j=1}^\infty (\eta_{n_j} f - 
\eta_{n_{j-1}}f) + \eta_1 f$ in the $L^p(h_\kappa^2)$ norm. Hence,
$$
  E_{2 n_j}(f)_{\kappa,p} \le \sum_{k = j+1}^\infty 
       \|\eta_{n_k}f - \eta_{n_{k-1}}f \|_{\kappa,p}. 
$$
Since $\eta_{n}$ preserves polynomials of degree $n$ and $\eta_n \eta_m f = 
\eta_m \eta_n f$ by definition, the triangle inequality shows that 
$$ 
\|\eta_{n_k}f - \eta_{n_{k-1}}f \|_{\kappa,p} \le 
  \|\eta_{2 n_k}f - \eta_{n_k} (\eta_{2 n_k}f)\|_{\kappa,p}
    + \|\eta_{2 n_k}f - \eta_{n_{k-1}} (\eta_{2 n_k} f) \|_{\kappa,p}.  
$$
Hence, using Proposition \ref{prop:3.7}, Theorem \ref{thm:3.2} and 
Proposition \ref{prop:3.11}, we get
\begin{align*}
\|\eta_{n_k}f - \eta_{n_{k-1}}f \|_{\kappa,p} 
   & \le c\, E_{n_{k-1}}(\eta_{2 n_k}f)_{\kappa,p}\\
   & \le c\, n_{k-1}^{-r} \|(- \Delta_{h,0})^{r/2} \eta_{2 n_k}f\|_{\kappa,p}\\
   & \le c\, n_k^{-r} \|(- \Delta_{h,0})^{r/2} \eta_{2 n_k}f\|_{\kappa,p}\\
   & \le c \,\omega_r(f;\pi/(2(2 n_k+\lambda)))_{\kappa,p}, 
\end{align*}
where the third inequality uses the fact that $n_k \le c n_{k-1}$. The 
definition of $n_k$ shows that 
$$
 \omega_r(f;\pi/(2(2 n_k+\lambda)))_{\kappa,p} \le  
    2^{j-k} \omega_r(f;\pi/(2(2 n_j+\lambda)))_{\kappa,p}, \qquad  k \ge j.   
$$
Consequently, we conclude that 
$$ 
E_{2 n_j}(f)_{\kappa,p} \le c \sum_{k = j+1}^\infty 
   \omega_r(f;\pi/(2(2 n_k+\lambda)))_{\kappa,p} 
 \le c \omega_r(f;\pi/(2(2 n_{j+1}+\lambda)))_{\kappa,p}. 
$$
Let $n \in \NN$. Choose a positive integer $j$ such that $n_j \le n/2 <
n_{j+1}$. Then it follows 
$$
  E_n(f)_{\kappa,p} \le c E_{2 n_j}(f)_{\kappa,p} 
   \le c \omega_r(f;\pi/(2(2 n_{j+1}+\lambda)))_{\kappa,p}  
   \le c \omega_r(f;\pi/(2(n +\lambda)))_{\kappa,p}  
$$
since $\omega_r(f;t)_{\kappa,p}$ is monotone nondecreasing. This proves
the proposition.  
\end{proof}

Now we can prove the equivalence of K-functional and modulus of smoothness:

\begin{thm}
For $f \in L^p(h_\kappa^2)$, $1 \le p \le \infty$, 
$$
c_1, \omega_r(f;t)_{\kappa,p} \le K_r(f;t)_{\kappa,p} \le 
 c_2\,  \omega_r(f;t)_{\kappa,p}.
$$
\end{thm}

\begin{proof}
The left side inequality is the Proposition \ref{prop:3.11}. To prove the 
right side, we choose $g = \eta_n f$ and use Proposition \ref{prop:3.12} to get
\begin{align*}
 K_r(f;t)_{\kappa,p} & \le \|f- \eta_n f\|_{\kappa,p} + 
  t^r \|(-\Delta_{h,0})^{r/2} \eta_n f\|_{\kappa,p} \\
 & \le E_n(f)_{\kappa,p} +  t^r \|(-\Delta_{h,0})^{r/2}\eta_n f\|_{\kappa,p}\\ 
 & \le c (1 + t^r n^{-r}) \omega_r(f; \pi/(2(n+\lambda)).  
\end{align*}
Taking $n = \inf  \{k \in \NN: \pi/(2(n+\lambda)) \le t\}$ and using the
monotonicity of $\omega_r(f;t)_{\kappa,p}$ proves the right side inequality.
\end{proof}

As a consequence of the above theorem, we are able to state:

\begin{cor} \label{cor:3.14}
For $f \in L^p(h_\kappa^2)$, $1 \le p \le \infty$, 
$$
 \omega_r(f;\delta t)_{\kappa,p} \le  c\, \max\{1,\delta^r\} 
    \omega_r(f;t)_{\kappa,p}.
$$
\end{cor}

This corollary allows us to replace the quantity $\pi /(2(n+\lambda))$ by
$n^{-1}$ in the Jackson type theorem. We state both the direct and the
inverse theorems.

\begin{thm}
For $f \in L^p(h_\kappa^2)$, $1 \le p \le \infty$, 
$$
   E_n(f)_{\kappa,p} \le c\, \omega_r(f;n^{-1})_{\kappa,p}.
$$
On the other hand, 
$$
 \omega_r(f,n^{-1})_{\kappa,p} \le c\, n^{-r} \sum_{k=0}^n (k+1)^{r-1} 
    E_k(f)_{\kappa,p}.
$$
\end{thm}

\begin{proof}
The proof of the direct theorem follows from Proposition \ref{prop:3.11}
and the above corollary. The inverse theorem follows from the Bernstein
type inequality using the standard argument. We note, however, even the 
inverse theorem needs the equivalence in Corollary \ref{cor:3.14}. The 
simple proof goes as follows: Let $P_n$ denote the polynomial of best 
approximation to $f$ of degree $n$. By the equivalent to the $K$-functional,
\begin{align*}
 \omega_r(f,n^{-1})_{\kappa,p} & \le c \left(\|f - P_{2^m}\|_{\kappa,p}
   + n^{-r} \|(-\Delta_{h,0})^{r/2} P_{2^m}\|_{\kappa,p} \right) \\
 & \le c\left( E_n(f)_{\kappa,p} + 
       n^{-r} \|(-\Delta_{h,0})^{r/2} P_{2^m}\|_{\kappa,p}\right).
\end{align*}
The rest of the proof follows from the standard argument of writting $P_{2^m} 
=  P_0 + \sum_{j=1}^m (P_{2^j} - P_{2^{j-1}})$ and using the Bernstein type 
inequality. 
\end{proof}

\begin{rem}
The proof of the direct theorem is not constructive because of the problem
that the equivalence of $\omega_r(f;t)_{\kappa,p}$ and $K_r(f;t)_{\kappa,p}$ 
is proved after the first direct estimate in Proposition \ref{prop:3.12}. 
However, at the end we see that the polynomial $\eta_n f$ satisfies the 
Jackson type estimate.
\end{rem}

\medskip

\begin{rem} \label{rem:3.2}
Let $\Lambda_\lambda$ be the differential operator defined by
\begin{equation} \label{eq:Lambda}
\Lambda_\lambda f(t)= w_\lambda^{-1}(t)\big[(1-t^2) w_\lambda(t) f'(t)\big]'.
\end{equation}
Then the Gegenbauer polynomials are eigenfunctions of $\Lambda_\lambda$; more 
precisely, the polynomials $C_n^\lambda$ satisfy the equation 
$\Lambda_\lambda f = -n(n+2 \lambda) f$ (\cite[p. 80]{Szego}). 
For $g \in L^1(w_\lambda, [-1,1])$, $g \sim \sum b_n C_n^\lambda$, we define
$(-\Lambda_\lambda)^{r/2}$ by
$$
 (-\Lambda_\lambda)^{r/2}g(t) \sim \sum b_n (n (n+2\lambda))^{r/2} 
  C_n^\lambda(t).
$$
Using Proposition \ref{prop:2.2} it is easy to verify that
$$
 (-\Delta_{h,0})^{r/2} P_n(h_\kappa^2;x,y) = V^{(y)}
    (-\Lambda_\lambda)^{r/2} P_n(w_\lambda;\langle x,y \rangle).
$$
So that fractional derivative of the reproducing kernel for $h$-harmonic
expansion is related to the fractional derivative of the kernel for Gegenbauer
expansion at the point $t=1$. See also Remark \ref{rem:4.1}.
\end{rem}

\section{Best approximation on the unit ball} 
\setcounter{equation}{0}

We consider weighted best approximation on $B^d$ for the weight function 
$$
 W_{\kappa,\mu}^B(x) = h_\kappa^2(x) (1-\|x\|^2)^{\mu-1/2}, \qquad x \in B^d,
$$
where $\mu \ge 1/2$ and $h_\kappa$ is an reflection invariant weight function 
defined on $\RR^d$. Let $a_{\kappa,\mu}$ denote the normalization constant for 
$W_{\kappa,\mu}^B$. Denote by $L^p(W_{\kappa,\mu}^B)$, $1 \le p \le \infty$, 
the space of measurable functions defined on $B^d$ with the finite norm
$$
  \|f\|_{W_{\kappa,\mu}^B,p}:= \Big(a_{\kappa,\mu} \int_{B^d} |f(x)|^p 
     W_{\kappa,\mu}^B(x) dx \Big)^{1/p}, \qquad 1 \le p < \infty,
$$
and for $p = \infty$ we assume that $L^\infty$ is replaced by $C(B^d)$, the 
space of continuous function on $B^d$. We consider 
$$
  E_n(f)_{W_{\kappa,\mu}^B,p} = \inf \{ \|f - P\|_{W_{\kappa,\mu}^B,p} : 
     P\in \Pi_n^d\}.   
$$
There is a close relation between best approximation on $B^d$ and on $S^d$. 

First we recall the relation between $h$-harmonics and orthogonal polynomials 
on the unit ball studied in \cite{X98a,X01a}. Let $\CV_n^d(W_{\kappa,\mu})$ 
denote the space of orthogonal polynomials of degree $n$ with respect to 
$W_{\kappa,\mu}$ on $B^d$. Elements of $\CV_n^d(W_{\kappa,\mu})$ are closely 
related to the $h$-harmonics associated with the weight function 
$$
h_{\kappa,\mu}(y_1,\ldots,y_{d+1}) = h_\kappa(y_1,\ldots,y_d)|y_{d+1}|^{\mu}
$$
on $\RR^{d+1}$, which is invariant under the group $G \times \ZZ_2$. Let 
$Y_n$ be such an $h$-harmonic polynomial of degree $n$ and assume that $Y_n$ 
is even in the $(d+1)$-th variable; that is, $Y_n(x,x_{d+1})=Y_n(x,-x_{d+1})$. 
We can write
\begin{equation} \label{eq:4.0}
   Y_n(y) = r^n P_n(x), \qquad  y = r(x,x_{d+1}) \in \RR^{d+1}, \quad 
 r = \|y\|, \quad (x,x_{d+1}) \in S^d,
\end{equation}
in polar coordinates. Then $P_n$ is an element of $\CV_n^d(W_{\kappa,\mu})$  
and this relation is an one-to-one correspondence \cite{X98a}. Furthermore, 
let $\Delta_h^{\kappa,\mu}$ denote the $h$-Laplacian associated with 
$h_{\kappa,\mu}$ and $\Delta_{h,0}^{\kappa,\mu}$ denote the corresponding
spherical $h$-Laplacian. When $\Delta_h^{\kappa,\mu}$ is applied to functions 
on $\RR^{d+1}$ that are even in the $(d+1)$-th variable, the spherical 
$h$-Laplacian can be written in polar coordinates $y = r(x,x_{d+1})$ as 
 (\cite{X01b}): 
$$
 \Delta_{h,0}^{\kappa,\mu} = \Delta_h - \langle x, \nabla 
 \rangle^2 - 2 \lambda \langle x, \nabla \rangle, \qquad 
    \lambda = \gamma_\kappa + \mu  + \frac{d-1}{2}. 
$$
in which the operators $\Delta_h$ and $\nabla=(\partial_1,\ldots,\partial_d)$ 
are all acting on $x$ variables and $\Delta_h$ is the $h$-Laplacian associated
with $h_\kappa$ on $\RR^d$. Define 
\begin{equation} \label{eq:4.1}
 D_{\kappa,\mu}^B := \Delta_h - \langle x, \nabla \rangle^2 - 
        2 \lambda \langle x, \nabla \rangle, 
\end{equation}
It follows that the elements of $\CV_n^d(W_{\kappa,\mu})$ are eigenfunctions 
of $D_{\kappa,\mu}^B$:
\begin{equation} \label{eq:4.2}
  D_{\kappa,\mu}^B P = - n (n+2\lambda) P,    
    \qquad P \in\CV_n^d(W_{\kappa,\mu}^B).
\end{equation}
For the classical weight function $W_\mu^B(x) = (1-\|x\|^2)^{\mu-1/2}$, the 
operator $D_{\kappa,\mu}^B$ becomes a pure differential operator which 
is classical (see \cite{AF} and \cite[Chapt. 12]{Er}). 

For $f \in L^2(W_{\kappa,\mu}^B)$, its orthogonal expansion is given by
$$
 L^2(W_{\kappa,\mu}^B) = \sum_{n=0}^\infty\bigoplus \CV_n^d(W_{\kappa,\mu}^B)
 : \qquad 
  f = \sum_{n=0}^\infty \proj_n^{\kappa,\mu} f, 
$$
where $\proj_n^{\kappa,\mu}:L^2(W_{\kappa,\mu})\mapsto\CV_n^d(W_{\kappa,\mu})$ 
is the projection operator. The fractional power of $D_{\kappa,\mu}^B$ on 
$f$ is defined by (see \eqref{eq:4.2})
$$
(-D_{\kappa,\mu}^B)^{r/2} f \sim  \sum_{n=0}^\infty (n (n+2\lambda))^{r/2}
   \proj_n^{\kappa,\mu} f, \quad  f \in L^p(W_{\kappa,\mu}^B).
$$
Using this operator we define
\begin{align*}
  \CW_r^p(W_{\kappa,\mu}^B):=\{f \in L^p(W_{\kappa,\mu}^B): 
     (-D_{\kappa,\mu}^B)^{r/2} \in L^p(W_{\kappa,\mu}^B)\}.
\end{align*}
Let $\|\cdot\|_{\kappa,\mu,p}$ denote the $L^p$ norm on $S^d$ with respect to
the weight function $h_{\kappa,\mu}$. We have the following important relation.

\begin{prop} \label{prop:4.1}
For $f \in \CW_r^p(W_{\kappa,\mu}^B)$, define $F(x,x_{d+1}) = f(x)$. Then 
$$
 \|(-D_{\kappa,\mu}^B)^{r/2} f\|_{W_{\kappa,\mu}^B,p} 
     = \|(- \Delta_{h,0}^{\kappa,\mu})^{r/2} F\|_{\kappa,\mu,p}. 
$$
\end{prop}

\begin{proof}
This follows from the equation \eqref{eq:4.1} and an elementary integral
$$
\int_{S^{d}}  g(y) d\omega(y) = \int_{B^d} \left[g(x, \sqrt{1-\|x\|^2})+
   g(x, -\sqrt{1-\|x\|^2})\right] \frac{d x}{\sqrt{1-\|x\|^2}}, 
$$
since the correspondence \eqref{eq:4.0} shows that 
$Y_n(h_{\kappa,\mu}^2; F) = \proj_n^{\kappa,\mu} f$. 
\end{proof}

In particular, this implies the following Bernstein type inequality:

\begin{cor} \label{cor:4.2}
For $P \in \Pi_n^d$ and $1 \le p \le \infty$, 
$$
 \|(-D_{\kappa,\mu}^B)^{r/2} P\|_{W_{\kappa,\mu}^B,p} \le c n^r 
  \|P\|_{W_{\kappa,\mu}^B,p}. 
$$
\end{cor}

\begin{proof}
Let $f = P$ in Proposition \ref{prop:4.1}. Then the Bernstein type inequality 
\eqref{eq:3.3} on $S^d$ shows 
$\|(-D_{\kappa,\mu}^B)^{r/2} P\|_{W_{\kappa,\mu}^B,p} \le c n^r 
    \|F\|_{\kappa,\mu,p}$,
where $F(x,x_{d+1}) = P(x)$. The fact that $\|F\|_{\kappa,\mu,p} = 
\|P\|_{W_{\kappa,\mu}^B,p}$ finishes the proof.  
\end{proof}

The K-functional between $L^p(W_{\kappa,\mu}^B)$ and 
$\CW_r^p(W_{\kappa,\mu}^B)$ is defined by 
$$
 K_r(f;t)_{W_{\kappa,\mu}^B,p}:= \inf \{ \|f-g\|_{W_{\kappa,\mu}^B,p} + 
  t^r \|(-D_{\kappa,\mu}^B)^{r/2} g\|_{W_{\kappa,\mu}^B,p} \},
$$
which can be used to characterize the best approximation by polynomials:

\begin{thm} \label{thm:4.3}
For $f \in L^p(W_{\kappa,\mu}^B)$, $1 \le p \le \infty$, 
$$
   E_n(f)_{W_{\kappa,\mu}^B,p} \le c\, K_r(f;n^{-1})_{W_{\kappa,\mu}^B,p}.
$$
On the other hand, 
$$
 K_r(f;n^{-1})_{W_{\kappa,\mu}^B,p} \le c\, n^{-r} \sum_{k=0}^n (k+1)^{r-1} 
    E_k(f)_{W_{\kappa,\mu}^B,p}.
$$
\end{thm} 

\begin{proof}
For $f \in L^p(W_{\kappa,\mu}^B)$, we associate it with a function $F$ on
$S^d$ defined by $F(x,x_{d+1}) = f(x)$, $(x,x_{d+1}) \in S^d$. Let $Y_n$ be a 
polynomial best approximation to $F$ of degree $n$ in $L^p(h_{\kappa,\mu}^2)$.
We can assume that $Y_n$ is even in its $(d+1)$-th variable. Indeed, let 
$\sigma_{d+1} Y_n(x,x_{d+1}): = Y_n(x,-x_{d+1})$; then  
$$
 \|F - (Y_n + \sigma_{d+1} Y_n)/2\|_{\kappa,\mu,p} \le 
   \|F - Y_n\|_{\kappa,\mu,p}/2 +\|F- \sigma_{d+1}Y_n \|_{\kappa,\mu,p}/2    
   = \|F - Y_n\|_{\kappa,\mu,p},
$$
where in the last step we changed the sign in the integral and used the fact 
that $F$ is even in its last variable. Thus, if $Y_n$ is a best approximation
to $F$, then so is $(Y_n + \sigma_{d+1} Y_n)/2$. Since $Y_n$ is even in its 
last variable, we can use $x_{d+1}^2 = 1 - \|x\|^2$ to write $Y_n (y)=r^n 
P_n(x)$, $y = r(x,x_{d+1})$, in which $P_n(x)$ is a polynomial of degree $n$
with $x \in B^d$. The integral in the proof of Proposition \ref{prop:4.1} 
shows then 
$$
\|f - P_n\|_{W_{\kappa,\mu}^B,p}
 = \|F-Y_n\|_{\kappa,\mu,p} = E_n(F)_{\kappa,\mu,p}.
$$ 
Hence, if $f \in \CW_r^p(W_{\kappa,\mu}^B)$, then by Theorem \ref{thm:3.2} and
Proposition \ref{eq:4.1},
$$
E_n (f)_{W_{\kappa,\mu}^B,p} \le  E_n(F)_{\kappa,\mu,p} \le c\, n^{-r} 
 \|(-\Delta_{h,0}^{\kappa,\mu})^{r/2}F
  \|_{\kappa,\mu,p} = c \,n^{-r}\|(-D_{\kappa,\mu}^B)^{r/2}f\|_{\kappa,p},  
$$
from which the direct estimate follows from the triangle inequality. The
inverse estimate again follows from the Bernstein type inequality. 
\end{proof}

\begin{cor}
For $f \in \CW_r^p(W_{\kappa,\mu}^B)$, $1 \le p \le \infty$, 
$$
E_n(f)_{W_{\kappa,\mu}^B,p} \le c\, n^{-r}
   \|(-D_{\kappa,\mu}^B)^{r/2}\|_{W_{\kappa,p}^B,p},   
$$
\end{cor}

This appears to be new even in the case of the classical weight function 
$W_\mu(x) = (1-\|x\|^2)^{\mu-1/2}$. Recall that $D_{\kappa,\mu}^B$ is a 
second order differential operator for $W_\mu$. It follows that if the 
$2s$ derivatives of $f$ are in $L^p(W_\mu)$, then the error of the best
approximation by polynomials of degree $n$ is in the order of $n^{-2s}$. 

One can define a mean on $B^d$ that corresponds to the spherical means 
$T_\theta^\kappa f$ on $S^d$, but the definition does not look natural 
on $B^d$. The same applies to the modulus of smoothness defined using 
$T_\theta^\kappa f$. One interesting question is to find another modulus of 
smoothness on $B^d$ that is in some sense natural and also equivalent to 
the K-functional defined above.

\medskip

\begin{rem} \label{rem:4.1}
In the case of $d=1$, $h_\kappa(x)=|x|^{2\kappa}$, $\kappa \ge 0$, and 
we get results on approximation theory with respect to the weight function 
$$
 w_{\kappa,\mu}(t) = |t|^\kappa (1-t^2)^{\mu -1/2}, \qquad -1< t<1.
$$
Even in this case the result in this section appears to be new. Only the
case $\kappa =0$, corresponding to the weight function for the Gegenbauer 
polynomials, has been studied in the literature. In particular, recall the
operator $\Lambda_\lambda$ defined in \eqref{eq:Lambda}, the Corollary 
\ref{cor:4.2} gives a Bernstein type inequality, which seems to be known 
only when $r$ is an even integer.

\begin{prop}
Let $1 \le p \le \infty$. For any polynomial $g$ of degree $n$ on $\RR$, 
$$
\|(-\Lambda_\lambda)^{r/2} g\|_{w_\lambda,p} \le c n^r \|g\|_{w_\lambda,p}. 
$$
\end{prop}
\end{rem}

\section{Best approximation on the simplex} 
\setcounter{equation}{0}

We consider weighted best approximation on the simplex $T^d$ for the weight 
function
$$
W_{\kappa,\mu}^T(x) = h_\kappa^2(\sqrt{x_1}, \ldots,\sqrt{x_d})
   (1-|x|)^{\mu-1/2} /\sqrt{x_1 \cdots x_d},  
$$
where $\mu \ge 1/2$ and $h_\kappa$ is an reflection invariant weight function 
defined on $\RR^d$ and $h_\kappa$ is even in each of its variables. The last 
requirement essentially limits the weight functions to the case of 
group $\ZZ_2^d$, for which  
\begin{equation}\label{eq:5.1}
 W_\kappa^T(x) = x_1^{\kappa_1-1/2} \cdots x_d^{\kappa_d-1/2}
   (1-|x|)^{\kappa_{d+1} -1/2} 
\end{equation}
(setting $\mu = \kappa_{d+1}$), which is the classical weight function on 
$T^d$, the case of hyperoctahedral group (see \eqref{eq:1.4}) and the case of 
$d =2$ and even dihedral group (see \eqref{eq:1.5}). The case of $h_\kappa$ in 
\eqref{eq:1.3} for the symmetric group, however, is excluded since it is not
even in its variables. 

The background on orthogonal expansion and approximation on $T^d$ is similar
to the case of the unit ball $B^d$. The definitions of various notions, such
as $\|\cdot\|_{W_{\kappa,\mu}^T,p}$, $E_n(f)_{W_{\kappa,\mu}^T,p}$ and 
$\CV_n(W_{\kappa,\mu}^T)$, are exactly the same as in the previous section
with $T^d$ in place of $B^d$. 

There is a close relation between orthogonal polynomials on $B^d$ and those 
on $T^d$ (\cite{X98a,X98b}). Let $P_{2n}$ be an element of $\CV_{2n}
(W_{\kappa,\mu}^B)$ and assume that $P_{2n}$ is even in each of its variables.
Then we can write $P_{2n}$ as $P_{2n} (x) = R_n(x_1^2,\ldots,x_d^2)$.
It turns our that $R_n$ is an element of $\CV_n(W_{\kappa,\mu}^T)$ and the 
relation is an one-to-one correspondence. In particular, applying 
$D_{W_{\kappa,\mu}^B}$ on $P_{2n}$ leads to a second order 
differential-difference operator acting on $R_n$. We denote this operator 
by $D_{\kappa,\mu}^T$. Then (\cite{X01b})
\begin{equation} \label{eq:5.2}
  D_{\kappa,\mu}^T P = - n (n+\lambda) P,    
    \qquad P \in\CV_n^d(W_{\kappa,\mu}^B), \quad \lambda = 
     \gamma_\kappa + \mu + \frac{d-1}{2}.
\end{equation}
For the  weight function \eqref{eq:5.1}, the operator is a second order 
differential operator, which takes the form 
$$
D_{\kappa,\mu}^T =
\sum_{i=1}^d x_i(1-x_i) \frac {\partial^2 P} {\partial x_i^2} - 
 2 \sum_{1 \le i < j \le d} x_i x_j \frac {\partial^2 P}{\partial x_i 
 \partial x_j} + \sum_{i=1}^d \left( \Big(\kappa_i +\frac{1}{2}\Big) -
    \lambda x_i \right) \frac {\partial P}{\partial x_i}    
$$
(recall $\mu = \kappa_{d+1}$ in this case). This is classical, already known
in \cite{AF} at least for $d=2$ (see also \cite[Chapt. 12]{Er}). For the 
formula of $D_{\kappa,\mu}^T$ in the case of hyperoctahedral group, see 
\cite{X01b}.

\begin{prop} \label{prop:5.1}
For $f \in L^p(W_{\kappa,\mu}^T)$, define $F(x) = f(x_1^2,\ldots,x_d^2)$. Then
$F \in L^p(W_{\kappa,\mu}^B)$ and  
\begin{equation} \label{eq:5.3}
\|(-D_{\kappa,\mu}^T)^{r/2} f\|_{W_{\kappa,\mu}^T,p} =  
     2^{-r} \|(-D_{\kappa,\mu}^B)^{r/2} F\|_{W_{\kappa,\mu}^B,p}
\end{equation} 
In particular, 
\begin{equation} \label{eq:5.4}
  K_r(F;t)_{W_{\kappa,\mu}^B,p} \le c \,K_r(f;t)_{W_{\kappa,\mu}^T,p}.
\end{equation}
\end{prop}

\begin{proof}
Writing $f = \sum c_k R_k$, where $R_k \in \CV_n^d(W_{\kappa,\mu}^T)$, we have
\begin{align*}
(-D_{\kappa,\mu}^T)^{r/2} f(x_1^2,\ldots,x_d^2) & = 
  \sum c_k (k(k+\lambda))^{r/2} R_k(x_1^2,\ldots,x_d^2) \\
 & = 
 2^{-r} \sum c_k (2k(2k+2\lambda))^{r/2} P_{2k}(x) \\&= 2^{-r}  
   (-D_{\kappa,\mu}^B)^{r/2} F(x),  
\end{align*}
where $P_{2k}(x) = R_k(x_1^2,\ldots,x_d^2) \in \CV_{2n}^d(W_{\kappa,\mu}^B)$. 
The elementary integral
$$
\int_{B^d} f(x_1^2,\ldots,x_d^2) dx = \int_{T^d} f(x) \frac{d x}
    {\sqrt{x_1 \cdots x_d}}  
$$
then proves equation \eqref{eq:5.3}. Furthermore, in the definition of 
$K_r(F;t)_{W_{\kappa,\mu}^B,p}$, the infimum is taken over $g \in 
\CW_r^p(W_{\kappa,\mu}^B)$. Taking infimum over functions that are even in 
each of its variables leads to an inequality, 
$$
 K_r(F;t)_{W_{\kappa,\mu}^B,p} \le \inf \{\|F-G\|_{W_{\kappa,\mu}^B,p} + 
   t^r \|(-D_{\kappa,\mu}^B)^{r/2} G\|_{W_{\kappa,\mu}^B,p}
$$
in which the infimum is taken over all $G \in \CW_r^p(W_{\kappa,\mu}^B)$ such
that $G(x) = g(x_1^2,\ldots,x_d^2)$. Since \eqref{eq:5.3} shows that 
$G\in \CW_r^p(W_{\kappa,\mu}^B)$ is equivalent to $g\in 
\CW_r^p(W_{\kappa,\mu}^T)$, the stated inequality follows.
\end{proof}

\begin{cor} \label{cor:5.2}
For $R \in \Pi_n^d$, $1 \le p \le \infty$,
$$
\|(-D_{\kappa,\mu}^T)^{r/2} R\|_{W_{\kappa,\mu}^T,p} \le c\,n^r  
   \|R\|_{W_{\kappa,\mu}^T,p}.
$$
\end{cor}

\begin{proof}
Let $P(x) = R(x_1^2,\ldots,x_d^2)$. The above proposition and 
Corollary \ref{cor:4.2} shows
\begin{align*}
 \|(-D_{\kappa,\mu}^T)^{r/2} R\|_{W_{\kappa,\mu}^T,p}
  & = 2^{-r} \|(-D_{\kappa,\mu}^B)^{r/2} P\|_{W_{\kappa,\mu}^B,p} \\
  & \le c\, (2n)^{r}\|P\|_{W_{\kappa,\mu}^B,p}  
   \le c\, n^{r}\|R\|_{W_{\kappa,\mu}^T,p},  
\end{align*}
which proves the stated Bernstein type inequality.
\end{proof}

\begin{thm} \label{thm:5.3}
For $f \in L^p(W_{\kappa,\mu}^T)$, $1 \le p \le \infty$, 
$$
  E_n(f)_{W_{\kappa,\mu}^T,p} \le c\, K_r(f;n^{-1})_{W_{\kappa,\mu}^T,p}.
$$
On the other hand, 
$$
 K_r(f,n^{-1})_{W_{\kappa,\mu}^T,p} \le c\, n^{-r} \sum_{k=0}^n (k+1)^{r-1} 
    E_k(f)_{W_{\kappa,\mu}^T,p}.
$$
\end{thm} 

\begin{proof}
For $f \in L^p(W_{\kappa,\mu}^T)$, we define $F(x) = f(x_1^2,\ldots,x_d^2)$. 
Then $F \in L^p(W_{\kappa,\mu}^B)$. Let $P \in \Pi_{2n}^d$ be a polynomial of
best approximation to $F$ in $L^p(W_{\kappa,\mu}^B)$. Let $P(x \varepsilon) 
= P(\varepsilon_1 x_1, \ldots, \varepsilon_d x_d)$ for $\varepsilon \in 
\{-1,1\}^d$. Since 
\begin{align*}
 \|F - 2^{-d} \sum_{\varepsilon \in \{-1,1\}^d}
  P((\cdot)\varepsilon) \|_{W_{\kappa,\mu}^B,p}
&\le 2^{-d} \sum_{\{-1,1\}^d} 
  \|F-P((\cdot)\varepsilon)\|_{W_{\kappa,\mu}^B,p}\\
& = \|F-P\|_{W_{\kappa,\mu}^B,p},
\end{align*}
where in the last step we changed the signs in the integral and used the fact 
that $F$ is even in each of its variables. Hence, we can assume that $P$ is 
even in each of its variables. Consequently, we can write $P(x) = 
R(x_1^2,\ldots,x_d^2)$ for a polynomial $R \in \Pi_n^d$. Hence, by 
Theorem \ref{thm:4.3}, 
$$
  E_n(f)_{{W_{\kappa,\mu}^T,p}} \le 
     \|f - R\|_{W_{\kappa,\mu}^T,p}
  = \|F-P\|_{W_{\kappa,\mu}^B,p} \le c K_r (F; (2n)^{-1})_{W_{\kappa,\mu}^B,p}.
$$ 
Thus, inequality \eqref{eq:5.3} proves the direct estimate. The inverse 
estimate again follows from the Bernstein type inequality. 
\end{proof}

\begin{cor}
For $f \in \CW_r^p(W_{\kappa,\mu}^T)$, $1 \le p \le \infty$, 
$$
E_n(f)_{W_{\kappa,\mu}^T,p} \le c\, n^{-r}
   \|(-D_{\kappa,\mu}^T)^{r/2}\|_{W_{\kappa,p}^T,p},   
$$
\end{cor}

\begin{rem} \label{rem:5.1}
In the case of $d=1$, $W_{\kappa,\mu}^T(x) = x^{\kappa-1/2}(1-x)^{\mu-1/2}$
is just the Jacobi weight function on $[0,1]$. The Jacobi weight function is
usually denoted by
$$
 w^{(\alpha,\beta)}(x) = (1-x)^{\alpha}(1+x)^{\beta}, \qquad -1 \le x \le 1.
$$
The result in this section appears to be new even in this case. In particular,
the operator $D_{\kappa,\mu}^T$ becomes the differential operator, denoted by 
$\Lambda_{\alpha,\beta}$, for the Jacobi polynomials
$$
 \Lambda_{\alpha,\beta} f(t) = \big[(1-t^2) w^{(\alpha,\beta)}(t) f'\big]'
      / w^{(\alpha,\beta)}(t).
$$
Then Corollary \ref{cor:5.2} gives a Bernstein type inequality, which seems 
to be known only when $r$ is an even integer.

\begin{prop}
Let $1 \le p \le \infty$. For any polynomial $g$ of degree $n$ on $\RR$, 
$$
\|(-\Lambda_{\alpha,\beta})^{r/2} g\|_{w^{(\alpha,\beta)},p} 
  \le c n^r \|g\|_{w^{(\alpha,\beta)},p}. 
$$
\end{prop}
\end{rem}

\bigskip\noindent
{\it Acknowledgment.} The author thanks a referee for his valuable comments
and for pointing out an oversight in Section 2.2 of the manuscript.


\begin{thebibliography}{99} 

\bibitem{AF} 
	P. Appell \and J. K. de F\'eriet.
	\textit{Fonctions hyperg\'eom\'etriques et hypersph\'eriques, 
        Polynomes d'Hermite}, Gauthier-Villars, Paris, 1926.

\bibitem{AW} 
        R. Askey and S. Wainger, 
        On the bexhavior of special classes of ultraspherical expansions,
        \textit{J. Anal. Math.}, \textbf{15} (1965), 193-244.

\bibitem{BBP} 
        H. Berens, P. L. Butzer and S. Pawelke, 
        Limitierungsverfahren von Reihen mehrdimensionaler Kugelfunktionen
        und deren Saturationsverhalten, 
        \textit{Publ. Res. Inst. Math. Sci. Ser. A.} \textbf{4} (1968), 
        201-268.

\bibitem{BL}
        H. Berens and Luoqing Li,         
        On the de la Vall\'ee Poussin means on the sphere,
        \textit{Results in Math.}, \textbf{24} (1993), 12-26.         

\bibitem{CZ} 
        A. P. Calderon and A. Zygmund
        On a problem of Mihlin,
        \textit{Trans. Amer. Math. Soc.}, \textbf{78} (1955), 209-224.  

\bibitem{DL} 
        R. DeVore and G. G. Lorentz,
        \textit{Constructive Approximation},  
        Springer-Verlag, Berlin, 1993

\bibitem{D1} 
        C. F. Dunkl,  
	Differential-difference operators associated to reflection groups,
        \textit{Trans. Amer. Math. Soc.} \textbf{311} (1989), 167--183.

\bibitem{D2} 
	C. Dunkl,
	Integral kernels with reflection group invariance,
 	\textit{Canad. J. Math.} \textbf{43} (1991), 1213-1227.

\bibitem{DX}
        C. F. Dunkl and Yuan Xu,
        \textit{Orthogonal polynomials of several variables},
        Cambridge Univ. Press, 2001. 
 
\bibitem{Er}
	A. Erd\'elyi, W. Magnus, F. Oberhettinger and F. G. Tricomi,
	\textit{Higher transcendental functions}, 
	McGraw-Hill, New York, 1953.

\bibitem{Kam}
        A. I. Kamzolov,
        The best approximation on the classes of functions $W_p^\alpha(S^n)$ by
        polynomials in spherical harmonics,
        \textit{Mat. Zametki}, \textbf{32} (1982), 285--293; English transl in 
        \textit{Math Notes}, \textbf{32} (1982), 622-628.

\bibitem{LW}
        Luoqing Li and Kunyang Wang,
        \textit{Harmonic analysis and approximation on the unit sphere}
        Science Press, Beijing, 2000.
 
\bibitem{LX}
        Zh.-K, Li and Yuan Xu,
        Summability of orthogonal expansions I, on unit sphere, and II, 
        on ball and simplex, \- {\it submitted}. 

\bibitem{LN}
        P. I. Lizorkin and S. M. Nikolskii,
        Approximation theory on the sphere,
        \textit{Proc. Steklov Inst. Math.}, \textbf{172} (1987),       
        295-302.         
 
\bibitem{P}
        S. Pawelke, 
        \"Uber Approximationsordnung bei Kugelfunktionen und algebraischen
        Polynomen, 
        \textit{T\^ohoku Math. J.}, \textbf{24} (1972), 473-486.

     
\bibitem{Ros} 
	M. R\"osler,
	Positivity of Dunkl's intertwining operator,
	\textit{Duke Math. J.\/},  \textbf{98} (1999), 445--463.

\bibitem{Rus}
        Kh. Rustamov, 
        On approximation of functions on the sphere,  
        \textit{Russian Acad. Sci. Izv. Math.}, \textbf{43} (1994), 311-329.

\bibitem{Rudin}
        W. Rudin,
        \textit{Real and Complex Analysis},
        McGraw-Hill, Inc., Boston, 1987.  

\bibitem{SW} 
 	E. M. Stein and G. Weiss,
 	\textit{Introduction to Fourier analysis on Euclidean spaces},
 	Princeton Univ. Press, Princeton, NJ, 1971.

\bibitem{Szego}
	G. Szeg\"{o},
	\textit{Orthogonal Polynomials},  
	Amer. Math. Soc. Colloq. Publ. Vol.23, Providence, 4th edition,
        1975.

\bibitem{V}
	N. J. Vilenkin,
	\textit{Special functions and the theory of group representations},
        American Mathematical Society Translation of Mathematics Monographs 
        \textbf{22}, 	
        American Mathematical Society, Providence, RI, 1968.

\bibitem{X97a}
	Yuan Xu,  
	Orthogonal polynomials for a family of product weight functions on the 
	spheres, \textit{Canad. J. Math.} \textbf{49} (1997), 175-192.  

\bibitem{X97b}
	Yuan Xu,
	Integration of the intertwining operator for $h$-harmonic polynomials
	associated to reflection groups, 
	\textit{Proc. Amer. Math. Soc.} \textbf{125} (1997), 2963--2973.

\bibitem{X98a}
	Yuan Xu,  
        Orthogonal polynomials and cubature formulae on spheres and on
        balls,  
        \textit{SIAM J. Math. Anal.} \textbf{29} (1998), 779--793.

\bibitem{X98b}
	Yuan Xu,  
        Orthogonal polynomials and cubature formulae on spheres and on 
        simplices,
        \textit{Methods Anal. and Appl.} \textbf{5} (1998), 169--184.

\bibitem{X00}
	Yuan Xu,            
        Funk-Hecke formula for orthogonal polynomials on spheres and on
        balls,
        \textit{Bull. London Math. Soc.} \textbf{32} (2000), 447--457.

\bibitem{X01a}
	Yuan Xu,            
	Orthogonal polynomials and summability in Fourier orthogonal 
        series on spheres and on balls,  
        \textit{Math. Proc. Cambridge Phil. Soc.}, \textbf{31} (2001), 
        139-155.    

\bibitem{X01b} 
	Yuan Xu,            
	Generalized classical orthogonal polynomials on the ball and on 
	the simplex,  
        \textit{Constr. Approx.}, \textbf{17} (2001), 383-412.

\bibitem{X02} 
	Yuan Xu,            
        Approximation by means of $h$-harmonic polynomials on the unit sphere,
        \textit{Adv. in Comp. Math}, to appear.

\end{thebibliography}
\enddocument